\documentclass[12pt]{amsart}
\usepackage{amssymb}
\usepackage{epsf}

\input epsf.tex

\def\myitemmargin{
\leftmargini=30pt 
\leftmarginii=21pt 
\leftmarginiii=19pt 
\leftmarginiv=19pt  
}

\newdimen\xsize
\newdimen\oldbaselineskip
\newdimen\oldlineskiplimit
\def\restorelineskip{\baselineskip=\oldbaselineskip%
\lineskiplimit=\oldlineskiplimit}
\def\putm[#1][#2]#3{
\hbox{\vbox to 0pt{\parindent=0pt%
\vskip#2\xsize\hbox to0pt{\hskip#1\xsize $#3$\hss}\vss}}}%
\def\putt[#1][#2]#3{
\vbox to 0pt{\noindent\hskip#1\xsize\lower#2\xsize%
\vtop{\restorelineskip#3}\vss}}

\DeclareFontFamily{U}{rsf}{\skewchar\font'177}%
\DeclareFontShape{U}{rsf}{m}{n}{<-6>rsfs5<6-8>rsfs7<8->rsfs10}{}%
\DeclareFontShape{U}{rsf}{b}{n}{<-6>rsfs5<6-8>rsfs7<8->rsfs10}{}%
\DeclareMathAlphabet\RSFS{U}{rsf}{m}{n}
\SetMathAlphabet\RSFS{bold}{U}{rsf}{b}{n}
\let\scr=\rfs
\def\mib#1{\boldsymbol{#1}} 
\def\sf#1{{\mathsf{#1}}} 

\def\slsf#1{{\slshape \sffamily #1\/}}

\def\msmall#1{\mathchoice{\hbox{\small$\displaystyle {#1}$}}{#1}{#1}{#1}}

\def\3{\ss}
\def\cc{{\mathbb C}}

\def\rr{{\mathbb R}}
\def\nn{{\mathbb N}}
\def\pp{{\mathbb P}}
  
\def\zz{{\mathbb Z}}

\def\aut{{\sf{Aut}}}

\def\codim{\sf{codim}\,}
\def\dscr{\sf{Dscr}}
\def\deg{\sf{deg}\,}

\def\dim{\sf{dim}\,}

\def\es{{\sf{es}}}

\def\length{\sf{length}\,}
\def\sym{\sf{Sym}}

\def\sfh{\sf{H}}

\def\hom{\sf{Hom}\,}
\def\scrhom{{\scrh\!\!o\!m}}

\def\ker{\sf{Ker}\,}
\def\lim{\mathop{\sf{lim}}}

\def\length{\sf{length}}

\def\mod{\;\sf{mod}\,}

\def\bfm{{\boldsymbol{m}}}
\def\max{\sf{max}}

\def\ord{\sf{ord}}

\def\pr{\sf{pr}}

\def\rank{\sf{rank}}
\def\res{\sf{Res}}

\def\Sing{\sf{Sing}}

\def\Sl{\mib{Sl}}
\def\pgl{\mib{PGl}}
\def\Gl{\mib{Gl}}

\def\bff{{\mib{F}}}

\def\bfn{{\mib{N}}}

\def\bfs{{\mib{S}}}
\def\bft{{\mib{T}}}

\def\eps{\varepsilon}
\def\epsi{\varepsilon}

\def\<{\langle}\let\la=\<
\def\>{\rangle}\let\ra=\>
 \let\bs=\bss

\def\d{\partial}

\def\ddef{\mathrel{{=}\raise0.3pt\hbox{:}}}
\def\deff{\mathrel{\raise0.3pt\hbox{\rm:}{=}}}
\def\inv{^{-1}}

\def\fraction#1/#2{\mathchoice{{\msmall{ #1\over#2}}}%
{{ #1\over #2 }}{{#1/#2}}{{#1/#2}}}

\def\leq{\leqslant}
\def\leq{\leqslant}

\def\xarr#1{\xrightarrow{#1}}
\def\emptyset{\varnothing}
\def\scirc{\mathop{\mathchoice{\hbox{\small$\circ$}}{\hbox{\small$\circ$}}%
{{\scriptscriptstyle\circ}}{{\scriptscriptstyle\circ}}}}
\def\longpoints{\leaders\hbox to 0.5em{\hss.\hss}\hfill \hskip0pt}
\def\stateskip{\smallskip}
\def\state#1. {\stateskip\noindent{\bf#1. }} 
\def\statep#1. {\stateskip\noindent{\bf#1 }} 

\def\proof{\state Proof. \,}
\def\proofp{\statep Proof. \,}
\def\Chi{\raise 2pt\hbox{$\chi$}}
\def\eg{\,{\sl{e.g.\/\ \hskip1pt plus1pt}}}
\def\ie{\,{\sl i.e.}}

\def\sli{{\sl i)\/} \,}         \def\slip{{\sl i)}}
\def\slii{{\sl i$\!$i)\/} \,}   \def\sliip{{\sl i$\!$i)}}
\def\sliii{{\sl i$\!$i$\!$i)\/} \,}    
\def\sliv{{\sl i$\!$v)\/} \,}   
\def\slv{{\sl v)\/} \,}         \def\slvp{{\sl v)}}

\def\mat{\sf{Mat}}

\def\det{\sf{det}}

\def\barr#1{\mskip1mu\overline{\mskip-1mu{#1}\mskip-1mu}\mskip1mu}
\def\Chi{\raise 2pt\hbox{$\chi$}}

\let\phI=\phi\let\phi=\varphi\let\varphi=\phI\def\vphi{\varphi}
  \def\frA{{\mathfrak{A}}}

\def\scrc{\scr{C}}
\def\scrd{\scr{D}}
\def\scrh{\scr{H}}
\def\scre{\scr{E}}

\def\scrm{\scr{M}}

\def\scro{{\scr O}}

\def\scrs{{\scr S}}
\def\scru{{\scr U}}

\def\scrw{{\scr W}}
\def\scrx{{\scr X}}
\def\scry{{\scr Y}}
\def\scrz{{\scr Z}}

\def\eps{\varepsilon}
\def\epsi{\varepsilon}
\def\bs{\backslash}

\def\d{\partial}

\def\1{{1\mkern-5mu{\rom l}}}
\def\Id{\sf{Id}}

\def\geq{\geqslant}\def\geq{\geqslant}
\def\inv{^{-1}}
\let\wh=\widehat
\let\wt=\widetilde
\def\fraction#1/#2{\mathchoice{{\msmall{ #1\over#2}}}%
{{ #1\over #2 }}{{#1/#2}}{{#1/#2}}}

\def\leq{\leqslant}

\def\emptyset{\varnothing}

\def\ti#1{{\tilde{#1}}}

\def\qed{\ \ \hfill\hbox to .1pt{}\hfill\hbox to .1pt{}\hfill $\square$\par}

\def\comment#1\endcomment{}

 \belowdisplayskip=\abovedisplayskip

\def\lineeqqno(#1){\hfill\llap{\vbox to 10pt%
{\vss\begin{align} \eqqno(#1)\end{align}\vss}}\vskip1pt}

\def\Matrix{\begin{matrix}}
\def\Endmatrix{\end{matrix}}
\def\Cases{\begin{cases}}
\def\Endcases{\end{cases}}

\advance\headheight 1.2pt

\def\newsection[#1]#2{\section{#2} \vskip-15pt\label{#1}\vskip15pt}
   \def\refsection#1{\slsf{Section \ref{#1}}}
\def\newsubsection[#1]#2{\subsection{#2}\label{#1}}

\newtheorem{thm}{Theorem}[section]
   \def\newthm#1{\begin{thm} \label{#1}}
   \def\refthm#1{\slsf{Theorem \ref{#1}$\,$}}
\newtheorem{lem}[thm]{Lemma}
   \def\newlemma#1{\begin{lem} \label{#1}}
   \def\lemma#1{\slsf{Lemma \ref{#1}}}
\newtheorem{prop}[thm]{Proposition}
   \def\newprop#1{\begin{prop} \label{#1}}
   \def\propo#1{\slsf{Proposition \ref{#1}}}
\newtheorem{corol}[thm]{Corollary}
   \def\newcorol#1{\begin{corol} \label{#1}}
   \def\refcorol#1{\slsf{Corollary \ref{#1}}}
\newtheorem{defi}{Definition}[section]
   \def\newdefi#1{\begin{defi} \label{#1}\sl }
   \def\refdefi#1{\slsf{Definition \ref{#1}}}

\def\eqqno(#1){\label{(#1)}}
\def\eqqref(#1){\eqref{(#1)}}

\usepackage{hyperref}

\input xy
\xyoption{all}

\numberwithin{equation}{section}

\begin{document}

\myitemmargin 
\baselineskip =16.0pt plus 2pt

\title[The Severi problem for Hirzebruch surfaces]%
{The Severi problem for Hirzebruch surfaces }
\author[V.~Shevchishin]{Vsevolod V.~Shevchishin}
\address{Mathematisches Institut\\
 Abteilung f{\"u}r reine Mathematik \\
 Albert-Ludwigs-Universit{\"a}t \\
 Eckerstra{\ss}e 1\\
 D-79104 Freiburg im Breisgau\\
Germany} \email{sewa@email.mathematik.uni-freiburg.de}
\dedicatory{} \subjclass{}
\begin{abstract} 
 We prove that the locus of irreducible nodal curves on a given Hirzebruch
 surface $\bff_k$ of given linear equivalency class and genus $g$ is
 irreducible.
\end{abstract}
\maketitle
\setcounter{tocdepth}{2}

\setcounter{section}{-1} \pagebreak[1]

\newsection[intro]{Introduction}

In the famous \slsf{Anhang F} of his book \slsf{``Vorlesungen {\"u}ber
 algebraische Geometrie''} \cite{Sev}, F.~Severi offered a proof of the
statement that the locus of irreducible plane curves of degree $d$ having
the prescribed number of nodes $\nu$ and no other singularities is
connected. However, his argument, which involved degenerating the curve into
$d$ lines, is not correct. The problem was attacked by several authors, see
review of Fulton \cite{Ful}, and the correct proof was given by Harris
\cite{Ha}, following original ideas of Severi.

In this paper we consider the Severi problem for \slsf{complex} Hirzebruch
surfaces.  Recall that the Hirzebruch surface $\bff_k$ of \slsf{index $k$}
($k\geq0$) is the fiberwise projectivization of the vector bundle $\scro \oplus
\scro(k)$ over $\pp^1$ and is equipped with the projection $\pr: \bff_k \to
\pp^1$. There exists non-singular rational curve $C_0 \subset \bff_k$ (resp., $C_\infty \subset
\bff_k$) of self-intersection $k$ (resp., $-k$), and every irreducible curve
$C$ on $\bff_k$ except $C_\infty$ is linearly equivalent to $d\cdot [C_0] + f\cdot [F]$
with non-negative $d$ and $f$, where $[F]$ is the linear equivalency class of
a fiber of $\pr: \bff_k \to \pp^1$. For integers $d \geq1$, $f\geq0$, and $0\leq g \leq
g_{\max}$ with $g_\max \deff \frac{k\,d(d-1)}{2} + (d-1)(f-1)$, let $\scrm^\circ=
\scrm^\circ(\bff_k, d,f, g)$ be the locus of irreducible nodal curves in $\bff_k$
in the linearly equivalency class $d\cdot [C_0] + f\cdot [F]$ of geometric genus $g$,
and $\scrm= \scrm(\bff_k, d,f, g )$ its closure. In the case $k=0$, when
$\bff_0 = \pp^1 \times \pp^1$, we additionally assume that $f\geq1$. It is easy to
show that $\scrm^\circ$ is a Zariski open subset of $\scrm$. The main purpose of
this paper is

\medskip
\state Theorem 1. {\it The variety $\scrm(\bff_k, d,f, g)$ is irreducible}.
\medskip

The main component of the proof is the following

\newthm{thm0.1} Each irreducible component of $\scrm(\bff_k, d,f, g)$ contains
a curve $C^\times$ which is a union of $d$ sections $C_i$, such that $C_i^2 =k$,
and $f$ pairwise distinct fibers of the ruling $\pr: \bff_k \to \pp^1$.
\end{thm}

The meaning of \refthm{thm0.1} is that there are no ``unexpected'' components
of the variety $\scrm(\bff_k, d,f, g)$, whereas each ``expected'' one can be
obtained by smoothing of an appropriate collection of nodes on $C^\times$. Using
the natural toric action on $\bff_k$, it is easy to show that each component
of $\scrm(\bff_k, d,f, g)$ contains a curve which consists of $d$ sections
$C_i$ as in \refthm{thm0.1} and a fiber $F$ with multiplicity $f$.  In this
way we are led to the main technical part of the paper which is the
\slsf{local Severi problem for ruled surfaces}. By this we mean the question
of description of possible nodal deformations of a given curve $C^*$ on a
ruled surface in a neighborhood of its unique compact component $F$ which is a
fiber of a ruling. The obtained solution allows to prove \refthm{thm0.1}
rather easily. Proving \slsf{Theorem 1}, we consider the action of the
monodromy group on the set of nodes of a curve $C^+$ obtained from the curve
$C^\times$ as above by smoothing some collection of nodes, such that $C^+ \in
\scrm(\bff_k, d,f, g=0)$.

\medskip 
The author's motivation for studying of the local Severi problem was its
relation to the \slsf{symplectic isotopy problem}. The techniques developed
in the papers \cite{Sh-1} and \cite{Sh-2} (see also \cite{Si-Ti}) allow to
show that every nodal pseudoholomorphic curve $C$ in $\bff_0 = \pp^1 \times \pp^1 \cong
S^2 \times S^2$ of genus $g\leq3$ is symplectically isotopic to an algebraic curve and
give evidences to hope that similar property holds for every pseudoholomorphic
curve $C$ in an arbitrary ruled surface $X$ provided $c_1(X) \cdot [C] >0$. So the
irreducibility of $\scrm(\bff_k, d,f, g)$ implies that the symplectic isotopy
class of such a curve is determined by the homology class and the genus of
$C$.  In particular, we have the following

\medskip 
\state Corollary. {\it There exists a unique symplectic isotopy class
 of irreducible nodal pseudoholomorphic curves in $S^2 \times S^2$ of given
 bi-degree $(d_1,d_2)$ and genus $g\leq3$.} 
\medskip

\medskip\noindent 
\slsf{Acknowledgments.} The author is strongly indebted to Urs~Hartl,
V.~Kharlamov, Vik.~Kulikov, and B.~Siebert for numerous valuable remarks and
suggestions which helped to clarify the problem and refine the exposition.
Many other valuable remarks and suggestions were made by H.~Flenner,
G.-M.~Greuel, Chr.~Lossen, St.~Nemirovski, St.~Orevkov, and E.~Shustin.



\baselineskip =14.0pt plus .5pt \nobreak


\vskip1cm

\tableofcontents

\newsection[sec:1]{Local Severi problem for ruled surfaces}

\newsubsection[1.1]{Moduli spaces of curves on ruled surfaces.}  Let us start
with a brief discussion of the working category for moduli spaces of curves.
First, we notice that the problem itself can be posed also in the case of the
ground field $\Bbbk$ of the non-zero characteristic, the answer could be quite
different, however. As an example of possible reasons, let us observe that in
the case $\mathsf{char} (\Bbbk)=2$ the discriminant of a polynomial of the
form $P(z,w)= a_0(z)w_0^2 + a_1(z) w_0w_1 + a_2(z) w_1^2$ with respect to $w$
is $a_1(z)^2 - 4 a_0(z)a_2(z) = a_1(z)^2$.  So in the contrast to the case
$\mathsf{char} (\Bbbk)=0$, every zero of the discriminant has multiplicity $2$.
This means that the method to distinguish the locus of nodal curves in a given
variety of curves used in \lemma{lem1.3} does not work in this case, at least
without appropriate changes.

This explains our restriction to the case of the field $\cc$ of
complex numbers as the ground field. So we can freely use all tools of the
complex analysis as the classics do \cite{Gr-Ha}. 

For a complex manifold $X$ and a complex curve $C$ with the smooth boundary,
we denote by $\scrh(C, X)$ the space of holomorphic maps $u:C \to X$ which
extend continuously up to the boundary $\d C$. In particular, $\scrh(C) \deff
\scrh(C, \cc)$ is the space of holomorphic functions which are continuous up
to boundary.

Denote by $\Delta$ the unit disc $\{ z \in \cc : |z| < 1\}$ and fix a
coordinate $w = [w_0 : w_1]$ on $\pp^1$. Denote by $\pr: \Delta \times \pp^1
\to \Delta$ the natural projection on the first factor.

\newdefi{def1.1} A \slsf{Weierstra\3 polynomial on $\Delta \times \pp^1$ of degree $d$}
is a polynomial of the form $P(z,w) = \sum_{i=0}^d a_i(z) w_0^{d-i} w_1^i$
whose coefficients $a_i(z)$ are holomorphic in $z\in \Delta$.  Its
\slsf{discriminant} with respect to $w$ is denoted by $\dscr(P)$.  $P(z,w)$
is \slsf{normalized} if $a_0(z)$ is a unital polynomial with zeroes in
$\Delta$.

A \slsf{curve in $\Delta \times \pp^1$ of degree $d$} is the zero divisor of some
Weierstra\3 polynomial $P(z,w) \not \equiv 0$. Such a $P$ is a \slsf{defining
 polynomial} for $C$. Observe that $C$ can be reducible and can have multiple
components; those could be only vertical lines $\ell_z \deff \{z\} \times \pp^1$; the
only possible compact components are also vertical lines $\ell_z$.

A curve is \slsf{nodal} if all its singular points are nodes.

A Weierstra\3 polynomial $P= \sum_{i=0}^d a_i(z) w_0^{d-i} w_1^i$ is
\slsf{proper} if its coefficient $a_i(z)$ lie in $\scrh(\Delta)$ and both
$a_0(z)$ and $\dscr(P)$ do not vanish on $\d \Delta$. The curve $C$ defined
by such $P(z,w)$ is also called \slsf{proper}.  The space of proper curves
of degree $d$ is denoted by $\scrz_d$.
\end{defi}

Two Weierstra\3 polynomials $P$ and $\ti P$ define the same curve $C$ iff
$\ti P = h \cdot P$ for some \slsf{invertible} $h(z) \in \scro(\Delta)$.
Thus every proper curve $C$ in $\Delta \times \pp^1$ can be represented by a
unique normalized Weierstra\3 polynomial, denoted by $P_C(z,w)$. 

\state Example. \slsf{Case $d=0$.} In this case a proper curve $C$ is given
by the Weierstra\3 polynomial $P(z,w) = a_0(z)$ for some $a_0(z) \in
\scrh(\Delta)$ with no zeroes on the boundary $\d \Delta$. So $C = \cup_i\;
m_i \cdot \ell_{z_i}$ where $z_i$ are the zeroes of $a_0(z)$ and $m_i$ their
multiplicities.

\newlemma{lem1.1} \sli The space $\scrz_d$ is an open set in the Banach
manifold of collections $(a_0(z), \ldots, a_d(z))$ where $a_1, \ldots\allowbreak a_d \in
\scrh(\Delta)$ and $a_0$ is a unital polynomial.

\slii The space $\scrw^d$ of proper Weierstra\3 polynomials of degree $d$ is
an open subset in the space $\Big\{ P= \sum_{i=0}^d a_i(z) w_0^{d-i} w_1^i :
(a_0(z), \ldots, a_d(z)) \in \big(\scrh(\Delta)\big)^{d+1} \Big\}$. The
natural map $F: \scrw^d \to \scrz_d$ associating to each polynomial its zero
divisor is a holomorphic surjection, and the kernel of differential $dF_P$
admits a closed complement at each $P \in \scrw^d$.
\end{lem}

\proof The first part is trivial. The second ones is obtained easily from the
following assertions:

\newlemma{lem1.1a}
\begin{itemize}
\item Every $f(z) \in \scrh(\Delta)$ with no zeroes on $\d\Delta$ admits a unique
 decomposition $f(z) = p(z) \cdot g(z)$ where $p(z)$ is a unital polynomial with
 zeroes in $\Delta$ and $g(z)$ is an invertible element in $\scrh(\Delta)$.
\item The set of such $f(z)$ is open in $\scrh(\Delta)$ and the decomposition map
 $f(z) \mapsto (p(z), g(z))$ is holomorphic. \qed
\end{itemize}
\end{lem}

\newdefi{def1.2} Let $\scrz^\circ_{d,\nu}$ be the locus of proper curves $C \in
\scrz_d$ which are nodal with exactly $\nu$ nodes and have no multiple
components.  Denote by $\scrz_{d,\nu}$ the closure of $\scrz^\circ_{d,\nu}$
in $\scrz_d$.
\end{defi}

Since the only compact curves in $\Delta \times \pp^1$ are fibers $\ell_z = \pr\inv(z)$ of
the projection $\pr: \Delta \times \pp^1 \to \Delta$, the group of holomorphic automorphisms of
$\Delta \times \pp^1$ is the semi-direct product of the group $\aut(\Delta) \cong \Sl(2,\rr)$ of
automorphisms of $\Delta$ and the group
\[
\pgl(2, \scro(\Delta)) =  \Gl(2, \scro(\Delta)) / \scro^*(\Delta)\cdot\Id.
\]
Namely, every fiber preserving automorphism $g$ of $\Delta \times \pp^1$ is
given by
\begin{equation}\eqqno(gl2-act1)
(z, [w_0: w_1]) \mapsto \big(z,\; [g_{00}(z)w_0 + g_{01}(z)w_1:\;
g_{10}(z)w_0 + g_{11}(z)w_1]\big)
\end{equation}
for some matrix $g = \begin{pmatrix}g_{00}(z) & g_{01}(z)\\ g_{10}(z) &
 g_{11}(z)\end{pmatrix}$ with holomorphic coefficients $g_{ij}(z) \in \scro(\Delta)$
with the non-vanishing determinant $\det(g) = g_{00}(z) g_{11}(z) -
g_{01}(z) g_{10}(z)$.  The action of $\pgl(2, \scro(\Delta))$ on $\Delta
\times \pp^1$ is induced by the action of $\Gl(2, \scro(\Delta))$ on the
space of Weierstra\3 polynomials given by
\begin{equation}\eqqno(gl2-act2)
P =\sum_{i=0}^d a_i(z)w_0^{d-i} w_1^i \mapsto P \scirc g \deff \sum_{i=0}^d
a_i(z) \big(g_{00}(z)w_0 +g_{01}(z)w_1\big)^{d-i} \big(g_{10}(z)w_0 +
g_{11}(z)w_1\big)^i
\end{equation}
The same formula defines the action of the algebra $\mat(2, \scro(\Delta))$ of
holomorphic $2\times2$-matrices $g$ on Weierstra\3 polynomials, such that
\begin{equation}\eqqno(act-dscr)
\dscr(P \scirc g)  = \det(g)^{2d-2} \cdot \dscr(P)
\end{equation}

\newdefi{def1.3} A \slsf{Banach analytic set of finite definition (BASFD)} is
a subset in a Banach manifold which locally is a zero set of a finite number
of holomorphic function.
\end{defi}

We refer to the book of Ramis \cite{Ra} ({\sl Chapitre II}, {{\S}{\S}}\,3 and 4)
for the main properties of such sets. The most important of them, nice to
have in mind, are:

\newprop{prop1.3a} \sli The germ of a Banach analytic set $\scry$ of finite
definition at any point $y\in \scry$ has finitely many irreducible
components $\scry_i$, each of them also being also a BASFD.

\slii Such an irreducible component $\scry_i$ admits locally a finite proper
branched analytic covering over a closed submanifold of finite codimension
in the ambient Banach manifold.

\sliii Let $\scry^*\subset \scry$ be the subset of those points where
$\scry$ is a Banach manifold. Then $\scry^*$ is open and dense in $\scry$,
the compliment $\scry \bs \scry^*$ is again a BASFD, and the irreducible
components $\scry_i$ are locally the closure of the connected components of
$\scry^*$.

\sliv The notion of \slsf{codimension} of a BASFD is well-defined and
well-behaving. In particular, $\codim_y (\scry' \cap \scry'') \leq \codim_y
(\scry')+ \codim_y (\scry'')$ for any $y\in \scry' \cap \scry''$. Besides, a
BASFD can not be represented as a finite (even countable) union of BASFD's
of higher codimension.

\slv Let $\scrz$ be a BASFD which is irreducible at $\zeta \in \scrz$, $\scry \ni \zeta$
its analytic subset of finite codimension such that each irreducible component
$\scry_i$ of $\scry$ at $\zeta$ has the same codimension $k$, $Z$ a finite
dimensional analytic set, $\Phi: Z \to \scrz$ an analytic map, and $z \in \Phi\inv(\zeta)$ a
point. Then every irreducible component $Y_i$ of the fiber $Y \deff
\Phi\inv(\scry)$ has codimension $\codim_z(Y_i \subset Z) \leq k$.
\end{prop}

\newdefi{def1.3.0} A property $\frak{A}$ holds for a \slsf{generic point $y$}
of a BASFD $\scry$ if\/ $\frak{A}$ holds for every $y \in \scry \bs \scrw$ for
some BASFD $\scrw \subset \scry$ such that $\scry \bs \scrw$ is dense in
$\scry$.
\end{defi}

\medskip
Let us turn back to the discussion about the category for varieties of curves.
Observe that for $d\geq4$ there exists no finite dimensional complete family of
deformations of a proper curve $C \in \scrz_d$. So, in contrary to the case of
an isolated singularity (see \eg \cite{Sh-2}), we can not avoid consideration
of infinite dimensional families. The properties listed in \propo{prop1.3a}
insure that the Zariski-like topology based on BASFD's allows to work as in
finite-dimensional case. Moreover, in the forthcoming proofs, one can replace
the spaces $\scrz_d(\Delta)$ by finite-dimensional subspaces in which the
coefficients $a_i(z)$ are polynomial of a fixed sufficiently high degree $N$.
The reason for such a possibility is that the definition of various varieties
and loci used in the proofs are given in terms of \slsf{polynomial relations between
jets $j^k_{z_1}a_i(z), \ldots ,j^k_{z_n}a_i(z)$ of the coefficients of Weierstra\3
polynomials of curves}, such that the number $n$ of jets and the degree $k$ are
given explicitly and can be estimated by the numerical invariants of the
problem. The same allows gives another one possible algebraic approach, in
which we let the coefficients $a_i(z)$ vary in a local ring $\scro_{Z, z_0}
^{\mathsf{alg}}$ of germs regular functions at a non-singular points $z_0$ on
an algebraic curve curve $Z$.  Geometrically this means that we consider germs
of curves on the ruled surface $Z \times \pp^1$ at the fiber $\{z_0\} \times \pp^1$.

\smallskip

\newsubsection[sec1:1a]{Varieties of nodal curves on $\Delta \times \pp^1$}

\newdefi{def1.3a} A \slsf{multiplicity pattern} of degree $d$ and length $l$
is a non-increasing sequence $\bfm = (m_1, \ldots, m_l)$ of positive
integers such that $d = |\bfm| \deff \sum_i m_i$. A multiplicity pattern
$\bfm' = (m_1', \ldots, m'_{l'})$ is a \slsf{degeneration} of a multiplicity
pattern $\bfm = (m_1, \ldots, m_l)$ if there exists a surjective map $\vphi:
\{1,\ldots,l\} \to \{1,\ldots,l'\}$ such that $m'_i = \sum_{\vphi(j)=i} m_j$
for every $i= 1,\ldots,l'$. In particular, $l' \leq l$ and $|\bfm'| =
|\bfm|$. Such a degeneration is \slsf{strict} if $l'< l$, or equivalently,
if $\bfm' \neq \bfm$.

A polynomial $P(z)$ \slsf{has zeros of multiplicity pattern $\bfm = (m_1,
  \ldots, m_l)$} if $P(z)= a_0 \cdot \prod_{i=1}^l (z-z_i)^{m_i}$ with
pairwisely distinct zeroes $z_i$. The locus of \slsf{unitary} polynomials of
a given multiplicity pattern $\bfm$ is denoted by $A_\bfm^\circ$ and its
closure by $A_\bfm$.

The pattern having $\nu$ of $2$'s and $d-2\nu$ of $1$'s is denoted by
$\bfm(d,\nu)$, and the corresponding locus $A_{\bfm(d,\nu)}$ by $A_{d,\nu}$.
\end{defi}

\newlemma{lem1.3.0} The closure $A_\bfm$ of the locus $A_\bfm^\circ$ is an affine
subset in the affine space of all unitary polynomials of degree $d \deff
|\bfm|$
\[\textstyle
\Big\{ P(z) = z^d + \sum_{i=1}^d a_i z^{d-i} : (a_1,\ldots,a_d) \in \cc^d
\Big\}
\]
of dimension $\dim A_\bfm= l =\length (\bfm)$. The complement $A_\bfm \bs
A_\bfm^\circ$ is the union of all $A_{\bfm'}$ over all strict degenerations
$\bfm'$ of $\bfm$.
\end{lem}

\proof Consider the Vi{\`e}te map $f: \cc^d \to \cc^d$ associating to each
$d$-tuple $(z_1,\ldots,z_d)$ the coefficients $(a_1,\ldots,a_d)$ of the
unitary polynomial $P(z) = z^d + \sum_{i=1}^d a_i z^{d-i} = \prod_{i=1}^d
(z-z_i)$. The map $f$ can be viewed as the quotient of $\cc^d$ with respect
to the action of the symmetric group $\sym_d$ permuting the coordinates
$(z_1,\ldots,z_d)$. In particular, $f$ is algebraic and proper. It remains
to notice that each $A_\bfm$ is the image of a linear subspace of $\cc^d$.
\qed

\newlemma{lem1.3} \sli  The locus $\scrz_{d,\nu}$ is a Banach analytic subset of
$\scrz_d$ of pure codimension $\nu$.

\slii The complement $\scrz_{d,\nu} \bs \scrz^\circ_{d,\nu}$ has codimension
 $1$ in $\scrz_{d,\nu}$. In particular, $\scrz^\circ_{d,\nu}$ is dense in $\scrz_d$.

\end{lem}

\proof Let $C \in \scrz^\circ_{d,\nu}$ be a proper nodal curve, $P_C= \sum a_i(z)
w_0^{d-i}w_1^i$ its normalized Weierstra\3 polynomial, and $\ell_{z_1}, \ldots ,
\ell_{z_k}$ be its compact components, $\ell_{z_j} = \{ z_j\} \times \pp^1$ with some $z_i
\in \Delta$. Then every $a_i(z)$ must be divisible by $p(z) \deff \prod_{j=1}^k (z-
z_j)$. An easy but important observation is that any curve $C'\in \scrz_d$ lying
sufficiently close to $C$ is also nodal and has at most $\nu$ nodes. Moreover,
if $C'$ has also $\nu$ nodes, then it normalization is diffeomorphic to the
normalization of $C$. Thus $C'$ must have the same number of compact
components $\ell_{z_1'}, \ldots , \ell_{z_k'}$, each $\ell_{z_j'}$ lying close to the
corresponding $\ell_{z_j}$. Thus $\scrz^\circ_{d,\nu}$ is the disjoint union of the
sets
\begin{align*}
\scrz^\circ_{d,\nu,k} \deff& \{ C \in \scrz^\circ_{d,\nu} : \text{$C$ has
exactly $k$ compact components} \}.
\\
\noalign{Put} \scrz_{d,\nu,k} \deff& \text{ the closure of
$\scrz^\circ_{d,\nu,k}$.}
\end{align*}

\smallskip
First we show that each $\scrz_{d,\nu,0}$ is a Banach analytic set of
codimension $\nu$ in $\scrz_d$.  The assertion, however, follows from
\cite{Sh-2}, \slsf{Lemma 2.13}. More precisely, the following statements
were proved:

{\leftmargini=25pt
\begin{itemize}
\it
\item
 Let $C \in \scrz_d$ be a nodal proper curve with no vertical component,
 $P_C(z,w)$ its Weierstra\3 polynomial, $\dscr(P_C)$ the discriminant, $N
 \deff \ord(\dscr(P_C))$ the order of vanishing of $\dscr(P_C)$ in $\Delta$, and
 $\scrd_C= z^N + \sum_{i=1}^N c_i z^{N-i}$ the unique unital polynomial of degree
 $N$, such that $\dscr(P_C) = \scrd_C \cdot h$ with some invertible $h\in \scrh(\Delta)$.
 Then $N =N(C)$ is constant on every connected component of $\scrz_d$, and the
 map $F_N$ associating to $C \in \scrz_d$ the coefficients $(c_1, \ldots, c_N) \in
 \cc^N$ of $\scrd_C$ is a holomorphic local submersion.
 
\item
 The restricted map $F_N : \scrz^\circ_{d,\nu,0} \to \cc^N$ takes value in $A_{N,\nu}$
 (see \refdefi{def1.3a}) and its image $F_N (\scrz^\circ_{d,\nu,0})$ form a dense
 set in the union of certain irreducible components of the set
 $F_N\inv(A_{N,\nu})$.
\end{itemize}
}

It follows that $F_N\inv(A_{N,\nu})$ is an analytic subset of codimension
$\nu$ in $\scrz_d$ and $\scrz_{d,\nu,0}$ is a locally finite union of some
its components. To determine which of the components form $\scrz_{d,\nu,0}$
we use the following observation. If the discriminant $\dscr(P_C)$ has the
multiplicity $2$ at $z_0\in \Delta$, then one of the following cases occurs:
\begin{enumerate}
\item $C$ has a single singular point on $\ell_{z_0}$ which is a node, all
 branches of $C$ meet $\ell_{z_0}$ transversely;
\item $C$ has a vertical inflection point at some $p\in \ell_{z_0}$, all
  remaining branches of $C$ meet $\ell_{z_0}$ transversally at pairwisely
  distinct points;
\item $\ell_{z_0}$ has a simple tangency with $C$ at two points and meets $C$
  transversally at remaining points;
\item $C$ has degree $d=2$; $\ell_{z_0}$ is a vertical component of $C$, and the
 remaining part $C' \deff C \bs \ell_{z_0}$ meets $\ell_{z_0}$ transversally at two
 distinct points.
\end{enumerate}
Thus the curve $C \in F_N\inv(A^\circ_{N,\nu}) \subset \scrz_d$ belongs to
$\scrz^\circ _{d,\nu,0}$ iff for each root $z_j$ of $\scrd_C = F_N(C)$ the
configuration at the line $\ell_{z_j}$ is as in the case $(1)$.

\medskip So it remains to treat the loci $\scrz^\circ_{d,\nu,k}$ with $k>0$.
For this purpose we observe that for every $C \in \scrz^\circ_{d,\nu,k}$
with vertical components $\ell_{z_1}, \ldots , \ell_{z_k}$ the normalized
Weierstra\3 polynomial of the curve $C' \deff C \bs \cup_{j=1}^k \ell_{z_j}$
is $P_C(z,w) / p(z)$ with $p= \prod_{i=1}^k (z-z_i)$. The nodality condition
implies that $C'$ meets each $\ell_{z_i}$ transversely at exactly $d$
points, so that $C' \in \scrz^\circ_{d, \nu -dk, 0}$. The space of unital
complex polynomials $p(z)$ of degree $k$ with zeros in $\Delta$ is naturally
identified with the space of divisors on $\Delta$ of degree $k$ which, in
turn, is the symmetric power $\sym^k\Delta$. Thus $\scrz^\circ_{d,\nu,k}$ is
naturally imbedded in $\scrz^\circ_{d, \nu -dk, 0} \times \sym^d\Delta$,
such that the complement parameterizes certain degenerate curves. The list
of possible degenerations is short:
\begin{itemize}
\item[(a)] some $z_1, \ldots ,z_k \in \Delta$ coincide; in this case $p(z)$
has multiple roots and the discriminant $\dscr(p(z))$ vanishes;
\item[(b)] some of $\ell_{z_j}$ is tangent to $C'$; in this case $p(z)$ has a
 common root with the discriminant $\dscr(P_{C'})$.
\end{itemize}

It follows that the complement of $\scrz^\circ_{d,\nu,k}$ in $\scrz^\circ
_{d, \nu -dk, 0} \times \sym^d\Delta$ is a Banach analytic set of
codimension $1$. Finally, we observe that there exists a natural holomorphic
map $G: \scrz_{d, \nu -dk, 0} \times \sym^d\Delta \to \scrz_d$ associating
to a curve $C' \in \scrz_{d, \nu -dk, 0}$ with the defining Weierstra\3
polynomial $P_{C'}(z,w)$ and a unital polynomial $p(z)$ of degree $k$ the
curve $C$ given by $p(z) \cdot P_{C'}(z,w)$. The map $G$ inverts the
decomposition $C = C' \cup \bigcup_{i=1}^k \ell_{z_i}$ of curves $C \in
\scrz^\circ_{d,\nu,k}$. Thus $G$ induces the isomorphism between $\scrz_{d,
\nu - d k, 0} \times \sym^d\Delta$ and $\scrz_{d,\nu,k}$.

The lemma follows.
\qed

\newdefi{def1.4} The \slsf{virtual nodal number $\delta= \delta(C)$} of a proper
curve $C \in \scrz_d$ is the maximum of those $\nu$ such that $C \in
\scrz_{d,\nu}$.

A \slsf{maximal nodal deformation of $C$} is a nodal proper curve $C'$ lying
on the component $\scry$ of $\scrz^\circ_{d,\delta(C)}$ whose closure $\barr
\scry$ contains $C$.
\end{defi}

Recall that the virtual nodal number $\delta(C,p)$ at an isolated singular
point $p$ of a curve $C$ on a smooth complex surface $X$ is defined as the
maximal number of nodes on a small holomorphic deformation of the germ of
$C$ at $p$. In the case when $C$ has no compact components such a maximal
nodal deformation $C'$ can be constructed as follows: Take the normalization
$u: \wt C \to C \subset X$ of $C$ and let $u'$ is generic holomorphic
perturbation of $u'$. Then $C' \deff u'(\wt C)$ is nodal and has exactly
$\delta(C,p)$ nodes near each singular point $p \in C$.

\newthm{thm1.5} Let $C^* \in \scrz_d$ be a proper curve, $\ell_{z_1}, \ldots ,\ell_{z_m}$
its vertical components, each taken with the appropriate multiplicity, and
$C^{\dag}$ the union of its non-compact components. Then
\begin{equation}\eqqno(1.4)
\delta(C^*) = d \cdot m  + \sum_{p \in \Sing(C^{\dag})} \delta(C^{\dag}, p)
\end{equation}
\end{thm}

\proof First we show that the r.h.s.\ of \eqqref(1.4) is realizable. Let $u:
\wt C^{\dag} \to \Delta \times \pp^1$ be the normalization of $C^{\dag}$.
The properness condition on $C^*$ implies that the restricted projection
$\pr: C^{\dag} \to \Delta$ is a non-ramified covering near the boundary $\d
\Delta$. Consequently, the boundary of $\wt C^{\dag}$ consists of smooth
circle and the normalization map $u$ extends continuously up to the boundary
$\d \wt C^{\dag}$. Perturbing holomorphically the $\pp^1$-component of $u$
we obtain a map $u': \wt C^{\dag} \to \Delta \times \pp^1$ whose image $C'
\deff u'(\wt C^{\dag})$ is a nodal curve with $\delta(C^{\dag}) \deff
\sum_{p \in
   \Sing(C^{\dag})} \delta(C^{\dag}, p)$ nodes. Now make a generic shift of each vertical
component $\ell_{z_i}$. Then the obtained lines $\ell_{z'_i}$ are pairwisely
disjoint and each of them meets $C'$ transversely in $d$ points. Thus $C''
\deff C' \cup \bigcup_{i=1}^m \ell_{z'_i}$ is nodal and has $d \cdot k + \delta(C^{\dag})$ nodes as
desired.

Obviously, the properness of $C''$ is equivalent to that of $C'$. The latter
property can be proved as follows. By the construction, the intersection of
$C'$ with each $\ell_z$ is proper and has index $d$. Thus the map $\phi: \Delta \to
\sym^d\pp^1$ given by $\phi: z \mapsto C' \cap \ell_z \subset \ell_z \cong \pp^1$ is well-defined and
holomorphic. Moreover, $\phi$ extends continuously up to the boundary $\d\Delta$.
There exists the natural isomorphism between the symmetric power $\sym^d\pp^1$
and the space $\pp^d$, such that the space of homogeneous polynomials
$\sum_{i=0}^d a_i w_0^{d-i}w_1^i$ of degree $d$ with complex coefficients $(a_0,\ldots
a_d) \in \cc^{d+1}$ is identified with the space $\sfh^0(\pp^d, \scr0(1))$. By
Grauert's theorem, $\phi$ can be lifted to a holomorphic map $\ti\phi: z\in \Delta \mapsto
(a_0(z), \ldots ,a_d(z))\in \cc^{d+1}$, also continuous up to the boundary $\d\Delta$. The
components $(a_0(z), \ldots ,a_d(z))$ are the coefficients of a defining
Weierstra\3 polynomial of $C'$.

\smallskip%
Showing that the r.h.s.\ of \eqqref(1.4) can not be exceeded, we start with
the observation that it is sufficient to consider the case when the curve
$C^*$ has single vertical component, say $\ell_0$ over the origin $0\in
\Delta$.  Moreover, we may additionally assume that each non-compact
component $C_i^*$ of $C^*$ is a disc and the projection $\pr: C_j^* \to
\Delta$ is ramified only over $0\in \Delta$.  Thus $C^*_j$ meets $\ell_0$ at
a single point $p_j$. Denote the number of non-compact components $C_j^*$ of
$C^*$ by $b$, the normalized Weierstra\3 polynomial of $C_j^*$ by
$P_{C^*_j}$, the resultant of $P_{C^*_i}$ and $P_{C^*_j}$ with respect to
$w$ by $\res(P_{C^*_i}, P_{C^*_j})$, the degree of $\pr: C^*_j \to \Delta$
by $d_j$, the multiplicity of $\ell_0$ in $C^*$ by $m$, the intersection
index of $C^*_i$ and $C^*_j$ by $\delta_{ij}$, and set $\delta_j \deff
\delta(C^*_j, p_j)$. We can additionally suppose that the discriminant
$\dscr(P_{C^*})$ vanishes only at the origin $0\in \Delta$. This implies
that different components $C^*_i$ and $C^*_j$ of $C^*$ can meet only at
$\ell_0$.

Our main idea is to relate the singularities of $C^*$ with zeroes of the
discriminant $\dscr(P_{C^*})$. First, we observe that total order of
vanishing of $\dscr(P_{C^*})$ on $\Delta$ remains constant under small
perturbations of $C^*$. Further, the decomposition $C^* = m \cdot \ell_0
\cup \bigcup_{j=1}^b C^*_j$ implies that the normalized Weierstra\3
polynomial of $C^*$ is $P_{C^*}= z^m \cdot \prod_{j=1}^b P_{C^*_j}$.
Consequently,
\[
\dscr(P_{C^*}) = z^{2m(d-1)} \cdot \prod_{j=1}^b \dscr(P_{C^*_j}) \cdot
\prod_{1\leq i<j \leq j} \Big(\res(P_{C^*_i}, P_{C^*_j})\Big)^2,
\]
so that for the order of vanishing of the discriminant $\dscr(P_{C^*})$ we
obtain
\[
\ord(\dscr(P_{C^*})) = 2m(d-1) + \sum_{j=1}^b \ord(\dscr(P_{C^*_j})) +
2\sum_{1\leq i<j\leq b} \ord(\res(P_{C^*_i}, P_{C^*_j})).
\]
Now let $u_j : \Delta \to C^*_j \subset \Delta \times \pp^1$ be
parameterizations of $C^*_j$. Make a small deformation $u'_j$ of each $u_j$
perturbing only the $\Delta$-component of $u_j$ and leaving the
$\pp^1$-component unchanged. Then for a generic choice of such $u'_j$ the
curves $C'_j \deff u'_j(\Delta)$ will be maximal nodal deformations of
corresponding $C_j^*$ with $\delta_j$ nodes and will meat each other
transversely at $\delta_{ij}$ points. Furthermore, each projection $\pr:
C'_j \to \Delta$ will have $d_j-1$ simple branchings. Thus we can conclude
that
\begin{align*}
\ord(\dscr(P_{C_j})) & = 2\delta_j + d_j -1 ,
\\
\ord(\res(P_{C_i}, P_{C_j})) &= \delta_{ij},
\end{align*}
Further, observe that by the definitions above $\sum_{j=1}^b (d_j -1) = d-b$
and the virtual nodal number of $C^{\dag}$ is $\delta^{\dag} \deff
\delta(C^{\dag}) = \sum_{j=1}^b \delta_j + \sum_{1\leq i<j \leq b}
\delta_{ij}$. So the r.h.s.\ of \eqqref(1.4) equals $\delta^{\dag} + m\,d$
and
\begin{equation}\eqqno(1.5)
\ord(\dscr(P_C)) = d-b + 2( \delta^{\dag} + m\,(d -1)) .
\end{equation}

Now let $C^\# \in \scrz_d$ be a nodal curve with $\delta = \delta(C^*)$
nodes lying sufficiently close to $C^*$. Since $C^*$ has $b$ non-compact
components which are discs, its boundary $\d C^*$ consists of $b$ circles.
The boundary $\d C^\#$ must have the same structure, so the number $b^\#$ of
the non-compact components must be at most $b$, $b^\# \leq b$. Denote by
$m^\#$ the number of vertical components of $C^\#$. Then $m^\# \leq m$.
Further, let $C^\natural$ be the union of non-compact component of $C^\#$.
Then $C^\natural$ has $\delta^\natural \deff \delta(C^\natural)= \delta -
d\cdot m^\#$ nodes. Applying the Riemann-Hurwitz formula to the projection
$\pr: C^\natural \to \Delta$ we see that it must have at least $d -b^\#$
ramification points, counted with multiplicities. Each ramification point of
$\pr: C^\natural \to \Delta$ makes the input $1$ in the degree of the
discriminant $\dscr(P_{C^\natural})$, whereas each node of $C^\natural$
gives $2$. So
\begin{align*}
\ord\big( \dscr(P_C) \big) &= \ord\big( \dscr(P_{C^\#}) \big) = 2  m^\#
(d-1) +  \ord\big( \dscr(P_{C^\natural}) \big) \geq
\\
& \geq 2  m^\# (d-1) +  2 ( \delta - d m^\# ) + d - b^\#.
\end{align*}
Comparing with \eqqref(1.4) and taking into account the inequalities $b^\#
\leq b$, $m^\# \leq m$, and $\delta \geq \delta^{\dag} + m d$ we conclude
that we must have the equality in all cases. Thus we obtain the relations
\[
b^\# = b \qquad \text{and} \qquad m^\# = m.
\]
in addition to the formula \eqqref(1.4). \qed

\newdefi{def1.6} Let $C\in \scrz_d$ be a proper curve, $\ell_{z_1}, \ldots, \ell_{z_m}$ its
vertical components taken with appropriate multiplicities, $C^\natural$ the
union of non-compact components, and $\wt C^\natural$ the normalization of
$C^\natural$. The \slsf{normalization of $C$} is the abstract union $\wt C
\deff \wt C^\natural \sqcup \bigsqcup_i \ell_{z_i}$ considered as an
abstract curve and equipped with the natural \slsf{normalization map} $u:
\wt C \to C \subset \Delta \times \pp^1$.

\end{defi}

\newcorol{cor1.6.0} Let $C^*\in \scrz_d$ be a proper curve with a Weierstra\3
polynomial $P_{C^*}$ and $\wt C^*$ its normalization. Then the total
vanishing order of the discriminant of $P_{C^*}$ is
\begin{equation}\eqqno(1.4a)
\ord(\dscr(P_{C^*})) = d - \Chi(\wt C^*) + 2 \delta(C^*)
\end{equation}

\end{corol}

\newcorol{cor1.6} Let $C^* \in \scrz_d$ be a proper curve, $C^\#$ its maximal
nodal deformation, and
\[
C^* = \bigcup_{i=1}^{m^*} \ell_{z^*_i} \cup \bigcup_{j=1}^{b^*} C^*_j \qquad
C^\# = \bigcup_{i=1}^{m^\#} \ell_{z^\#_i} \cup \bigcup_{j=1}^{b^\#} C^\#_j
\]
their decomposition into irreducible components, such that $C^*_j$ and
$C^\#_j$ are non-compact ones. Then $b^* = b^\#$, $m^* = m^\#$, and,
possibly after a re-indexation, $C^*_j$ and $C^\#_j$ have the same geometric
genus.
\end{corol}

In other words, maximal nodal deformations of a given proper curve $C^*$ are
exactly those nodal curves which can be obtained by the following
construction: Each component is deformed preserving its geometric genus, in
particular, each compact component $\ell_z$ of $C^*$ is shifted in the
$z$-direction.

\newcorol{cor1.7} For a given proper curve $C^* \in \scrz_d$ with the nodal
number $\delta^* \deff \delta(C^*)$, the space $\scrz_{d, \delta^*}$ is
irreducible at $C^*$.

In particular, any two maximal nodal sufficiently small deformations $C',
C''$ of $C^*$ can be connected by a holomorphic family $C_\lambda \in
\scrz_{d,
 \delta^*}^\circ$, $\lambda \in \Delta$, also lying sufficiently close to
$C^*$.
\end{corol}

\proof Since every maximal nodal deformation is given by generic
``independent'' deformations of individual components of $C^*$, it is
sufficient to prove the special case when $C^*$ is irreducible.  The
subcases $d=0$ (in which $C^*$ is a vertical line $\ell_z$) is trivial. The
remaining cases $d\geq1$ follow from \cite{Sh-2}, \slsf{Lemma 1.9\,\,c)}.
\qed

\newdefi{def1.8} Let $C^*$ be a curve in $\Delta \times \pp^1$ and $\ell_z$ a
line, $z \in \Delta$. The \slsf{virtual nodal number $\delta(C^*, \ell_z)$
of
 $C$ at the line $\ell_z$} is the virtual nodal number of the restriction $C^*
\cap \big( \Delta(z, \varepsilon) \times \pp^1\big)$ of $C$ to a
sufficiently small neighborhood of $\ell_z$.
\end{defi}

\newcorol{cor1.8} For a given proper curve $C^* \in \scrz_d$, the virtual nodal
number $\delta(C^*, \ell_z)$ of $C$ at any line $\ell_z$ is well-defined and
\begin{equation}\eqqno(1.4b)
\delta(C^*) = \sum_{z\in \Delta} \delta(C^*, \ell_z).
\end{equation}

\end{corol}
\bigskip

\newsubsection[1.1a]{Equisingular families of curves}

\newdefi{def1a.1} Let $C^*\in \scrz_d$ be a proper curve and $z^* \in \Delta$ a
point. The space of \slsf{equisingular deformations of $C^*$ at the line
$\ell_{z^*}$} is the {\it connected} component $\scrz^\es_d(C^*, z^*)$ of the
locus $\scry^*$ of curves $C \in \scrz_d$ such that
\begin{itemize}
\item[(1)] the multiplicity of $\ell_{z^*}$ in $C$ and in $C^*$ coincide;
\item[(2)] the discriminants of the Weierstra\3 polynomials of $C^*$ and
 $C$ have the same zero order at $z^*$.
\end{itemize}
\end{defi}

\newprop{prop1a.2} The locus $\scrz^\es_d(C^*, z^*)$ is a Banach analytic set
of finite codimension in $\scrz_d$ and is irreducible at $C^*$.
\end{prop}

\proof Let $m^*$ be the multiplicity of $\ell_{z^*}$ in $C^*$, $a^*_0(z),
\ldots, a^*_d(z)$ the coefficients of the normalized Weierstra\3 polynomial
of $C^*$, and $n$ the order of vanishing of $\dscr(P_{C^*})$ at $z^*$. Then
$\scry^*$ is given by equations of vanishing of the jets $j^{n-1}_{z^*}
\dscr(P_{C^*})$ and $j^{m^*-1}_{z^*}a^*_i(z)$ with $i= 0,\ldots ,d$.  Thus
$\scry^* \subset \scrz_d$ is a BASFD. Since $\scrz^\es_d(C^*, z^*)$ is locally a
union of a finite number of components of $\scry^*$, it is also a BASFD.

\medskip
Dividing the Weierstra\3 polynomial $P_C$ of any curve $C \in \scry^*$ by
$(z- z^*)^{m^*}$, we reduce the problem of irreducibility of $\scry^*$ at
$C^*$ to the case $m^*=0$.

\smallskip
Now consider a holomorphic family of curves $C_\lambda$ in $\scry^*$, $\lambda \in \Delta$, such
that $C_0 = C^*$. Let $P_\lambda$ be the corresponding holomorphic family of the
normalized Weierstra\3 polynomials. Then for $|\lambda| \leq \epsi \ll 1$ the zero
divisor of $\dscr(P_\lambda)$ in $\Delta(z^*, \epsi)$ is $(z-z^*)^n$. Consequently, the
projections $\pr: C_\lambda \to \Delta$ are not ramified over the punctured disc $\check
\Delta(z^*, \epsi) \deff \Delta(z^*, \epsi) \bs \{z^*\}$. Moreover, the \emph{topological
 structure} of singularities of $C_\lambda$ at $\ell_{z^*}$ is constant in $\lambda$. The
later means that all topological (\ie, numerical) invariant describing the
structure of $C_\lambda$ at $\ell_{z^*}$ and their projections $\pr: C_\lambda \to \Delta$ coincide.
For example, $C_\lambda$ have the same number $l$ of local irreducible components at
$\ell_{z^*}$, say $C_1(\lambda), \ldots, C_l(\lambda)$, the same ramification degree $m_i$ of
projections $\pr: C_i(\lambda) \to \Delta$ over $z^*$, the same intersection indices
$\delta_{ij} = C_i(\lambda) \cap C_j(\lambda)$ at $\ell_{z^*}$, and so on. The constancy of these
numerical invariants follows from the fact, that otherwise for some $C_\lambda$
close to $C^*$ we would obtain a zero point $z'$ of some $\dscr(P_{C_\lambda})$,
which is close to but distinct from $z^*$.

\smallskip
Proceeding forth, we now observe the jet $j^{n -1}_{z^*} \dscr(P_C)$ is a
polynomial function of the $(n-1)$-jets $j^{n -1}_{z^*} a_i(z)$ of the
coefficients $a_i(z)$ of the Weierstra\3 polynomial $P_C$ at the point
$z^*$. This allows to reduce the irreducibility of $\scry^*$ at $C^*$ to the
following problem.

Let $Z_{d,n}$ be the space of Weierstra\3 polynomials $Q(z,w) = \sum_{i=0}^d
b_i(z)w_0^{d-i} w_1^i$ such that each $b_i(z)$ is a polynomial of degree at
most $n-1$. Define the projection $j^{n-1}: \scrz_d \to Z_{d,n}$ associating
to each normalized Weierstra\3 polynomial $P_C= \sum_{i=0}^d a_i(z)w_0^{d-i}
w_1^i$ its jet $j^{n-1}P_C(z,w) \deff \sum_{i=0}^d j^{n-1}_{z^*} a_i(z)\;
w_0^{d-i} w_1^i$. Set $Q^* \deff j^{n-1}P_{C^*}$ and define the set $Y^*
\subset Z_{d,n}$ by the equation $j^{n-1}_{z^*} \dscr Q =0$. Then $j^{n-1}:
\scrz_d \to Z_{d,n}$ is a holomorphic surjection, $Y^*$ is algebraic in
$Z_{d,n}$ and $\scry^* = \big(j^{n-1}\big)\inv(Y^*)$. Thus the
irreducibility of $\scry^*$ at $C^*$ and that of $Y^*$ at $Q^*$ are
equivalent.

The irreducibility of $Y^*$ would follow from the existence of a dominant
algebraic map $G: W \to Y^*$ with an irreducible variety $W$. Constructing
such a variety $W$, we first consider the space $V$ of polynomial maps
\[
f: t\in \cc \mapsto (f_1(t), \ldots , f_l(t)) \in \big( \pp^1 \big)^l
\]
such that each component $f_i(t): \cc \to \pp^1$ a fixed sufficiently large
degree $N$. For the ramification degrees $m_1, \ldots, m_l$ introduced
above, we set $\ti f_i(t) \deff (z^*+ t^{m_i}, f_i(t)) \in \cc \times \pp^1$
and consider the curves $\wt C_f \deff \cup_i \ti f_i(\cc) \cup m^*
\ell_{z^*} \subset \cc \times \pp^1$. If $\sum_{i=1}^l m_i =d$, then each
$C_f$ is a proper curve of degree $d$ in $\cc \times \pp^1$ and can be given
by a Weierstra\3 polynomial $P_f(z, w) = \sum_{i=0}^d a_i(z) w_0^{d-i}
w_1^i$ with polynomial $a_i(z)$ of degree $D = m^* + l \cdot N$. Define the
map $G: V \to Z_{d,n}$ associating to a polynomial map $f: \cc \to
(\pp^1)^l$ as above the $n-1$-jet $j^{n-1}P_f$ of the polynomial $P_f$. Then
$G: W \to Z_{d,n}$ is algebraic by the construction.

We contend that the image $G(V) \subset Z_{d,n}$ contains $Y$ and the
preimage $W \deff G\inv(Y)$ is irreducible, provided the degree $N$ of maps
$f_i: \cc \to \pp^1$ is chosen large enough. To show the first assertion,
let us consider a curve $C \in \scry^*$ close to $C^*$ and local irreducible
non-vertical components $C_i$ of $C$ at $\ell_{z^*}$.  Then each $C_i$
admits a unique local parameterization $\ti g_i: t \in \Delta(0, \epsi)
\mapsto \ti g_i(t) \in \Delta \times \pp^1$ of the form $\ti g_i(t) = (z^* +
t^{m_i}, g_i(t))$ with $g_i(t) \in \scro(\Delta(0, \epsi))$. Taking the
$N$-jet $f_i(t) \deff j^N_0g_i(t)$ we obtain a map $F: \scry^* \to W$ with
$F: C \mapsto (j^N_0g_1(t), \ldots j^N_0g_l)$ which is well-defined locally
near $C^* \in \scry^*$.  Moreover, it follows from the construction that $F$
is holomorphic and the composition $G \scirc F: \scry^* \to Z_{d,n}$
coincides with $j^{n-1} : \scry^* \to Z_{d,n}$ for $N$ large enough.

So it remains to show irreducibility of $V = G\inv(Y^*) \subset W$. The
crucial observation is that $V$ is given by \emph{linear} conditions on the
coefficients $c_{ij}$ of the components $f_i(t) = \sum_{j=0}^N c_{ij} t^j$
of $f=(f_1(t), \ldots f_l(t)) \in W$. Indeed, we have either the condition
$j^{d_{ij}}_0f_i(t) = j^{d_{ij}}_0f_j(t)$ of coincidence of the jets of
different $f_i(t)$ and $f_j(t)$ up to certain degree $d_{ij}$, or the
condition of vanishing of certain coefficients $c_{ij}$ of the components
$f_i(t) = \sum_{j=0}^N c_{ij} t^j$. The lemma follows.
\qed

\newlemma{lem1a.3} Let $C^* \in \scrz_d$ be a proper curve, $\delta^* \deff
\delta(C^*)$ its virtual nodal number, $m^*$ the multiplicity of $\ell_0$ in
$C^*$, and $b_0$ the number of local irreducible non-vertical components of
$C^*$ at $\ell_0$. Then the codimension $\codim\big( \scrz^\es_d(C^*, 0)
\subset \scrz_d \big)$ is at least $\delta^* + m^* + (d -b_0)+ 2$ except the
following cases:
\begin{enumerate}
\item[\sli] the non-vertical components of $C^*$ at $\ell_0$ are non-singular and
 disjoint from each other; in this case the codimension is $\delta^* + m^* + (d
 -b_0)$;
\item[\slii] the non-vertical components of $C^*$ at $\ell_0$ are non-singular
 and all of them except two are disjoint from each other; in this case the
 codimension is $\delta^* + m^* + (d -b_0) +1$;
\item[\sliii] the non-vertical components of $C^*$ at $\ell_0$ are disjoint from each
 other and all of them except one are non-singular; the horizontal projection
 of the exceptional component on the $Ow$-axis $\ell_0$ has degree $2$; in this case
 the codimension is $\delta^* + m^* + (d -b_0)+1$.
\end{enumerate}

Moreover, in the case \slii the local non-vertical branches $C_1^*$ and
$C_2^*$ of $C^*$ which meet $\ell_0$ at the same point $p \in \ell_0$ satisfy the
following condition: either $C_1^*$ and $C_2^*$ meet transversally at $p$ or
both $C_1^*$ and $C_2^*$ are tangent to $\ell_0$ at $p$.
\end{lem}

\proof Using the fact that $\scrz_{d, \delta^*}$ is irreducible, we obtain
\[
\codim\big( \scrz^\es_d(C^*, 0) \subset \scrz_d \big) = \codim\big(
\scrz^\es_d(C^*, 0) \subset \scrz_{d, \delta^*} \big) + \codim\big(
\scrz_{d, \delta^*} \subset \scrz_d \big).
\]
The last summand equals $\delta^*$ by \lemma{1.3}. The codimension
$\scrz^\es_d (C^*, 0) \subset \scrz_{d, \delta^*}$ can be estimated using
the following facts. First, according to the description of maximal nodal
deformations from \refcorol{cor1.6}, we must impose $m^*$ complex conditions
to obtain the multiplicity $m^*$ of $\ell_0$ in $C^*$. Second, let us denote
by $d_i$ the degrees of non-vertical local irreducible components $C^*_i$
($i=1,\ldots,b_0$) of $C^*$ at $\ell_0$. Then the ramification degree of the
projection $\pr: C^*_i \to \Delta$ is $d_i -1$. Since these ramifications
can ``walk'' in an arbitrary way under deformation of $C^*$ in $\scrz_{d,
\delta^*}$, we obtain additional $\sum_{i=1}^{b_0} d_i -1 = d - b_0$
conditions.

Compute the remaining parameters describing the locus equisingular
deformations.  Let $z=\zeta_i ^{d_i},\; w= \phi_i(\zeta_i)$ be a local
parameterization of the component $C^*_i$, $i=1,\ldots b_0$, $N$ a large
enough integer and $\psi_i(\zeta_i)$ polynomials of degree at most $N$ with
sufficiently small coefficients. Denote by $\psi$ the whole collection
$(\psi_1, \ldots,\psi_{b_0})$, by $C_{\psi,i}$ the curve with the local
parameterization $z=\zeta_i ^{d_i},\; w= \phi_i(\zeta_i)+ \psi_i(\zeta_i)$,
and by $C_\psi$ the curve $\cup_{i=1}^{b_0} C_{\psi,i} \bigcup m^* \ell_0$.
Then $C_\psi$ lie in $\scrz(\Delta(\eps))$ for some $\eps>0$ sufficiently
small. Furthermore, there exists a Weierstra\3 polynomial $P(z,w)$ of degree
$d$ which close enough to $P_{C^*}(z,w)$ and has the following properties:
\begin{itemize}
\item The $N$-jets of the coefficients of $P(z,w)$ at the origin $0\in \Delta$
 coincide with the corresponding jets of the Weierstra\3 polynomial $P_{C_\psi}$;
\item The $N$-jets of the coefficients of $P(z,w)$ and of $P_{C^*}(z,w)$
 coincide at every zero $z_j \neq 0$ of the discriminant $\dscr(P_{C^*})$.
\end{itemize}
By the construction, the curve given by such a polynomial $P(z,w)$ lies in
$\scrz_{d,\delta^*}$ have the same behavior at $\ell_0$ as $C_\psi$.

Let $p_1,\ldots p_l$ be the intersection points of $\ell_0$ with the local
non-vertical components $C^*_i$ of $C^*$ and $\mu_1, \ldots,\mu_l$ the
number of such components $C^*_i$ passing through $p_j$. The above
construction allows to move the components $C^*_i$ in the $w$-direction
separately. This gives $\sum_{j=1}^l (\mu_j -1)$ additional conditions
defining $\scrz^\es_d(C^*, 0)$ inside $\scrz_{d,
 \delta^*}$.

Let some local non-vertical components $C^*_i$ of $C^*$ at $\ell_0$ admits a
local parameterization $z=\zeta_i ^{d_i},\; w= \phi_i(\zeta_i)$ with $d_i
\geq 2$ and $\ord_{\zeta_i=0}(\phi_i(\zeta_i) - \phi_i(0)) \ddef s_i \geq
2$. Then for every non-zero polynomial $\psi(\zeta_i)$ of the form
$\psi(\zeta_i) = \sum_{j=1}^{s_i-1} c_j \zeta_i^j$ with sufficiently small
coefficients $c_j$ the curve with the parameterization $z=\zeta_i ^{d_i},\;
w= \phi_i(\zeta_i) + \psi(\zeta_i)$ will have singular points lying outside
the line $\ell_0$. Thus we obtain $s_i-1$ more parameter(s).

Since we are interested only in the case $\sum_{j=1}^l (\mu_j -1) + \sum
_{d_i \geq 2}  (s_i -1) =1$ we obtain exactly one of the possibilities \slii
and \sliii of the lemma.

Finally, assume that $C^*$ is as in the case \slii and that two local branches
$C^*_1$ and $C^*_2$ of $C^*$ passing through the same point $p$ on $\ell_0$ are
both non-vertical at $p$. Then $C^*_1$ and $C^*_2$ have local parameterizations
$w = \phi_i(z)$ with some holomorphic functions $\phi_1(z)$ and $\phi_2(z)$ which are
defined in a neighborhood of the origin $0\in\Delta$ and satisfy condition $\phi_1(0)
= \phi_2(0)$. Then the tangency condition of $C^*_1$ and $C^*_2$ at $p$ is given by
$\phi_1'(0) = \phi_2'(0)$. This is a complex condition which increases the
codimension of $\scrz_d^\es(C^*, 0)$ in $\scrz_{d,\delta^*}$ by 1. 
\qed

\newlemma{lem1.7a} Let $\scry$ be an irreducible BASFD and $\phi: \scry \to
\scrz_d$ a holomorphic map, such that for generic $y \in \scry$ the curve
$C_y
\deff \phi(y)$ is reducible. Then locally at every given $y^* \in \scry$ there exists
a holomorphic map $\Phi: \scry \to \prod_{j=0}^l \scry_j$ such that
\begin{enumerate}
\item $\scry_j$ is a local irreducible component of some space $\scrz_{d_j,
 \nu_j}$ at some $C^*_j$;
\item for every $y \in \scry$ with $\Phi(y) = (C_0, \ldots, C_l)$ and $C_y \deff \phi(y)$
 one has the decomposition $C_y = \cup_{j=0}^l C_j$; in particular, the curve
 $C^* \deff \phi(y^*)$ is decomposed into components $(C_0^*, \ldots, C_l^*)$ of
 $\Phi(y^*)$;
\item for a generic $y \in \scry$, $\Phi_0(y)$ is the union of vertical
 components of $C_y$ whereas the remaining $\Phi_j(y)$ ($j =1,\ldots,l$) are
 non-vertical irreducible components of $C_y$;
\item near the given $y^* \in \scry$ the map $\phi: \scry \to \scrz_d$ factorizes
 into the composition of $\Phi: \scry \to \prod_{j=0}^l \scry_j$ with the natural map
 $\prod_{j=0}^l \scry_j \subset \prod_{j=0}^l \scrz_{d_j, \nu_j} \to \scrz_d$ given by $(C_0,
 C_1, \ldots C_l) \mapsto \cup_{j=0}^l C_j$.
\end{enumerate}

Moreover, if\/ $\scry$ is a local irreducible component of $\scrz_{d, \nu}$
at some curve $C^*$, then the map $\Phi: \scry \to \prod_{j=0}^l \scry_j$ is
an isomorphism at $C^*$.
\end{lem}

\proof The map $\phi: \scry \to \scrz_d$ is given by a ``universal''
Weierstra\3 polynomial $P(z,w; y) = \sum_{i=0}^d a_i(z;y) w_0^{d-i} w_1^i$
whose coefficients $a_i(z,y)$ depend holomorphically on $z\in \Delta$ and
$y \in \scry$. The equation $P(z,w;y) =0$ defines the ``universal curve''
$\scrc \subset \Delta \times \pp^1 \times \scry$, whose fiber $\scrc_y \deff
\scrc \cap (\Delta \times \pp^1 \times \{ y\} )$ is the curve $C_y \deff
\phi(y)$.  It follows from the construction that $\scrc$ is a BASFD.

Take a point $p$ lying on a curve $C_{y_0} = \phi({y_0})$ corresponding to a
generic $y_0 \in \scry$ and consider local irreducible components of
$\scrc_i$ at points $p$.  The genericity of $y_0$ implies that each
$\scrc_i$ defines one local irreducible component of $C_{y_0}$ at $p$, and
also one local irreducible component on $C_{y'}$ for $y \in \scry$ close
enough to ${y_0}$. Consequently, locally in a neighborhood of $C_{y_0}$, the
irreducible components of $\scrc$ are correspond to the irreducible
components of $C_{y_0}$. If $y$ varies in sufficiently small neighborhood of
$y^*$, each non-compact component of $C_y$ stays close to exactly one
component of $C_{y^*}$. Thus the monodromy can interchange only compact
components of $C_{y_0}$. Hence for every non-compact component of $C_{y_0}$
we obtain one component $\scrc_i$ in a neighborhood of $C_{y^*}$. We define
$\scrc_0$ as the union of remaining components of $\scrc$. The case $\scrc_0
= \emptyset$ can occur and may be treated in the obvious way. By our
construction, for a generic $y \in \scry$ close to $y^*$ the intersection
$\scrc_0 \cap C_y$ is the union of all compact components of $C_y$.

For each component $\scrc_i$, let $d_i$ be the degree and $\nu_i$ the
virtual nodal number of the intersection $C_{y,i} \deff \scrc_i \cap C_y$
for generic $y$. Then $C_{y,i}$ lies in $\scrz_{d_i, \nu_i}$ and the induced
map $\Phi_i: \scry \to \scrz_{d_i, \nu_i}$. We set $\scry_i$ to be an
irreducible component of $\scrz_{d_i, \nu_i}$ containing $\Phi_i(\scry)$.
This construction extends---with the full accordance with definitions---also
to the case of $\scrc_0$ and gives the following. $d_0 = \nu_0 = 0$,
$\scrz_{0,0}$ is the set of all unitary polynomials $a_0(z)$ with zeroes in
$\Delta$, and $\scry_0$ the component of $\scrz_{0,0}$ containing
polynomials $a_0(z)$ of degree $m$ equal to the number of vertical
components of $C_{y,i}$ with generic $y\in \scry$.

The map $\Phi$ is given by its components $\Phi_0, \ldots ,\Phi_l$. The last
assertion of the lemma could be now seen easily. \qed

\newlemma{lem1a.4} Let $C^*,C^\#  \in \scrz_d$ be proper curves with the normalized
Weierstra\3 polynomials $P_{C^*}$ and $P_{C^\#}$, respectively. Assume that
\begin{itemize}
\item[(1)] $C^\#$ is close enough to $C^*$;
\item[(2)] both discriminants $\dscr(P_{C^*})$ and $\dscr(P_{C^\#})$ have zero only
 at the origin $0 \in \Delta;$
\item[(3)] the multiplicity $m^\#$ of the line $\ell_0$ in $C^\#$ is strictly less than
 the multiplicity $m^*$ of the line $\ell_0$ in $C^*$.
\end{itemize}

Then a maximal nodal deformation of $C^\#$ can be obtained from a maximal
nodal deformation of $C^*$ by smoothing appropriate nodes lying on vertical
lines.
\end{lem}

\proof The assertion, as the lemma itself, is trivial in the case $d=1$.
Thus we assume that $d\geq2$. Using \lemma{lem1.7a} we can reduce the
assertion of the lemma to the special case when $C^\#$ is irreducible.
Denote by $\delta^* \deff \delta(C^*)$ and $\delta^\# \deff \delta(C^\#)$
the corresponding virtual nodal numbers and by $m^*$ the multiplicity of
$\ell_0$ in $C^*$. Let $\scry$ be the irreducible component of $\scrz_{d,
\delta^\#}$ passing through $C^\#$. Then by property (2) and
\refcorol{cor1.6} every curve $C$ in $\scry$ has one non-vertical component.
Further, by condition $(1)$ implies that $\scry$ passes through $C^*$.

Consider the locus $\scry^*$ of those $C \in \scry$ for which the
discriminant $\dscr(P_C)$ vanishes only at the origin $0 \in \Delta$. Use
notation $y$ for an element in $\scry^*$ and $C_y$ for the corresponding
curve. Let $\hat\phi: \scry^* \to \scrz_d$ be the holomorphic map
associating to the normalized Weierstra\3 polynomial $P(z,w;y)$ of the curve
$C_y$ the Weierstra\3 polynomial $P(\zeta^d, w; y)\ddef \wh P(\zeta, w;y)$.
The geometric meaning of $\hat\phi$ is that $\wh P(\zeta, w;y)$ is the the
normalized Weierstra\3 polynomial of the curve which is the pre-image of
$C_y$ with respect to the map $F: \Delta \times \pp^1 \to \Delta \times
\pp^1$ given by $(\zeta, w) \in \Delta \times \pp^1 \overset{F}\longmapsto
(\zeta^d, w) \in \Delta \times \pp^1$. In particular, the preimage of $C^\#$
consists of $d$ discs every of which has degree $1$ in $ \Delta \times
\pp^1$. Application of \lemma{lem1.7a} yields $d$ holomorphic maps
$\hat\phi_i: \scry^* \to \scrz_1$, $i=1,\ldots,d$, such that $\hat\phi(y)$
is the union of curves $\hat\phi_i(y)$. On the level of Weierstra\3
polynomials we obtain $\wh P(\zeta, w;y)= \prod_{i=1}^d \wh P_i(\zeta,
w;y)$.

An easy---but crucial for us---observation is that the image of each curve
$\hat \phi_i(y)$, $i=1,\ldots,d$, with respect to the map $F$ is the curve
$C_y$ itself. The geometric meaning of this fact is as follows. For generic
$y \in \scry^*$ each $\hat\phi_i(y)$ is the graph of a map $w= f_i(\zeta)$
such that the curve $C_y$ admits the parameterization $z= \zeta^d,
w=f_i(\zeta)$. For arbitrary $y \in \scry^*$ the curve $\hat\phi_i(y)$ is
the union of the graph of some holomorphic map $w= f_i(\zeta)$ with a
vertical line $\ell_0$, taken with the multiplicity $m(y)$ equal to the
multiplicity of $\ell_0$ in $C_y$.

\medskip
Now consider the space $Y$ consisting of pairs $(q(\zeta), Q(w, \zeta))$
where
\begin{itemize}
\item $q(\zeta)$ is a unital polynomial of degree $d$ with zeroes in $\Delta$;
\item $Q(\zeta, w)= b_0(\zeta) w_0 + b_1(\zeta)w_1$ is a Weierstra\3 polynomial
of degree $1$ whose coefficients $b_0(\zeta),\; b_1(\zeta)$ are polynomials
of sufficiently high degree $N$, such that $b_0(\zeta)$ is unital.
\end{itemize}
For every such $\eta=(q,Q) \in Y$, we denote by $C_\eta$ the curve in
$\Delta \times \pp^1$ given by equations $z= q(\zeta)$ and $Q(\zeta, w) =0$.
Then $C_\eta$ depends algebraically on $\eta \in Y$, so that we obtain a
holomorphic map from $\Phi: Y \to \scrz_d$. Since $Y$ is irreducible, the
claim of the lemma follows from the fact that for $N$ large enough the
family $\{C_\eta\}_{\eta \in Y}$ contains sufficiently small equisingular
deformations of both $C^*$ and $C^\#$, as also their maximal nodal
deformations. In terms of the map $\hat\phi_1$ above, the first part of the
assertion means that the non-compact component of $C^*$ (resp. $C^\#$) can
be approximated by a curve given by the parameterization $z= \zeta^d,
w=f^*(\zeta)$ (resp. $z= \zeta^d, w=f^\#(\zeta)$) where $f^*(\zeta)$ and
$f^\#(\zeta)$ are polynomials of degree $N-m^*$ and $N$, respectively. The
existence of such a simultaneous approximation is evident. Moreover, we may
still assume that $C^\#$ is close to enough to $C^*$. This means that
corresponding points $\eta^*,\;\eta^\# \in Y$ are close to each other.
Finally, observe that a maximal nodal deformation of a given curve $C_\eta$
with $\eta =(q(\zeta), Q(\zeta, w)) \in Y$ close to $\eta^*$ can be obtained
by the following construction: A perturbation $\ti q(\zeta)$ of $q(\zeta)$
such that $\ti q(\zeta)$ has only simple branchings, and a generic
perturbation $\wt Q(\zeta, w)$ of $Q(\zeta, w)$ such that the curve given by
the equation $\wt Q(\zeta, w)=0$ is a maximal nodal deformation of the curve
given by $Q(\zeta, w)=0$.
\qed

\break

\newsubsection[1.3]{Local Severi problem for ruled surfaces}

\nobreak

\newthm{thm1.9} Let $C^* \in \scrz_d$ be a proper curve and $\delta^* \deff
\delta(C^*)$ its virtual nodal number. Then every local irreducible
component $\scry$ of $\scrz_{d, \nu}$ at $C^*$ ($\nu \leq \delta^*$)
contains $\scrz_{d, \delta^*}$.
\end{thm}

The meaning of the theorem is as follows. Fix any maximal nodal sufficiently
small deformation $C^{\dag} \in \scrz^\circ_{d, \delta^*}$ of $C^*$. Then
every irreducible component of $\scrz_{d, \nu}$ at $C^{\dag}$ ($\nu \leq
\delta^*$) can be reached by smoothing an appropriate collection of
$\delta^* -\nu$ nodes on $C^{\dag}$. So the theorem ensures that every
component of $\scrz_{d, \nu}$ at $C^*$ can also be obtained in this way.  In
particular, there are at most $\binom{\delta^*}{\nu}$ irreducible components
of $\scrz_{d, \nu}$ at $C^*$. Another interpretation is that any non-maximal
nodal sufficiently small deformation $C$ of $C^*$ can be degenerated into a
nodal curve with exactly one additional node.

\medskip Before giving the complete proof of \refthm{thm1.9} we consider
certain special cases.

\slsf{ Case 1: $\nu = \delta(C^*)$.} Here the claim of the theorem is
covered by the definition and \refcorol{cor1.7}. Thus we may assume that
$\nu < \delta(C^*)$.

\slsf{ Case 2: $d=1$.} In this case a curve $C^*$ is the zero set of the
polynomial $a_0^*(z) w_0 + a_1^*(z) w_1$ with a unital polynomial $a_0^*(z)$
and a holomorphic function $a_1^*(z) \in \scrh(\Delta)$. Let $\ell_{z^*_i}$
be the vertical components of $C^*$ and $m_i^*$ their multiplicities. Then
$\prod_i (z-z_i^*) ^{m_i^*}$ is the normalized greatest common divisor of
$a_0^*(z)$ and $a_1^*(z)$. Moreover, $\delta(C^*) = \sum_i m_i^*$ and the
curve is nodal iff every $m_i^*=1$. Further, a normalized Weierstra\3
polynomial of a sufficiently small deformation $C$ of $C^*$ is given by
small perturbations $a_0(z)$ and $a_1(z)$ of the coefficients $a_0^*(z)$ and
$a_1^*(z)$. So in the case of nodal $C$ with $\nu = \delta(C) < \delta(C^*)$
nodes we can deform the coefficients $a_0(z)$ and $a_1(z)$ in the way
preserving $\nu$ existing common zeroes of $a_0(z)$ and $a_1(z)$ and
creating $\delta(C^*) -\nu$ new ones.

\medskip
\slsf{ Case 3: A generic curve $C$ in $\scry$ is reducible.} Then
\lemma{lem1.7a} allows us to reduce \refthm{thm1.9} to the case when a
generic curve $C$ in $\scry$ is irreducible. Indeed, \lemma{lem1.7a}
provides a decomposition $C = \cup_{j=0}^l C_j$ of curves $C$ in $\scry$,
and a maximal nodal deformation of $C^*$ is the union of generic maximal
nodal deformation of individual pieces $C^*_j$ in the decomposition $C^* =
\cup_{j=0}^l C_j^*$.

\medskip
\slsf{ Case 4: The discriminant $\dscr(P_{C^*})$ has at least two distinct
 zeroes.}  Let $z^*_1, \ldots, z^*_l$ be the zero pairwise distinct points of the
discriminant $\dscr(P_{C^*})$, $l \geq 2$. By induction, we may assume that the
assertion of the theorem holds for all curve $C$ for which the total zero
order of the discriminant $\ord(\dscr(P_C))$ is strictly less than that for
$C^*$. In particular, it is so for any restriction of $C^*$ to a sufficiently
small neighborhood of any line $\ell_{z^*_i}$, \ie, for curves $C^* \cap
\big(\Delta(z^*_i, \varepsilon) \times \pp^1\big)$ with $\varepsilon$ small enough. Let us fix such a small
$\varepsilon$ and denote by $\Delta_i$ the disc $\Delta(z^*_i, \varepsilon)$. Further, fix a sufficiently
small neighborhood $\scru$ of $C^*$ in $\scrz_d(\Delta)$ and denote by $R_i : \scru
\to \scrz_d(\Delta_i)$ the restriction map associating to each curve $C$ its ``$i$-th
slice'' $C \cap \big(\Delta_i \times \pp^1\big)$. Fix a generic nodal curve $C^\circ \in \scry^\circ\cap
\scru$ lying sufficiently close to $C^*$. Then the curves $R_i(C^\circ)$ are also
nodal and the corresponding nodal number $\nu_i \deff \delta(R_i(C^\circ))$ are
independent of the choice of such a curve $C^\circ$. Moreover, $\sum_{i=1}^l \nu_i =\nu$.
Denote by $\scry_i$ the component of $\scrz_{d, \nu_i} (\Delta_i)$ which contains
$R_i(C^\circ)$.

Take a maximal nodal deformation $C^+$ of $C^*$ lying sufficiently close to
$C^*$. Then $R_i(C^+)$ is a maximal nodal deformation of $R_i(C^*)$. By the
inductive assumption, $R_i(C^+)$ belongs to $\scry_i$. This means that
smoothing certain collection of nodes on $C^+$ we obtain a curve $C'$ such
that $R_i(C')$ lie in $\scrz^\circ _{d, \nu_i} (\Delta_i)$. The theorem
follows from the next

\newlemma{lem1.10a} In the notation introduced above, let $C^\circ, C' \in \scrz_d$
are two nodal curves lying sufficiently close to $C^*$. Assume that for each
$i=1, ..., l$ the slices $R_i(C^\circ)$, $R_i(C')$ lie in the same component
$\scry_i$ of $\scrz_d(\Delta_i)$. Then $C^\circ$ and $C'$ lie in the same component
$\scry$ of $\scrz_d$.
\end{lem}

\proof We give the proof only for the case when $C^\circ$ and $C'$ have no
vertical components. The general case follows easily from this special one.

Consider a function $\bff(\vec \zeta, \vec p, G(z))$ which associates with
given pairwise distinct points $\zeta_1, \ldots, \zeta_n \in \Delta$, polynomials $p_1(z), \ldots,
p_n(z)$ of degrees $\deg(p_i(z)) = m_i-1$, respectively, and with a given
holomorphic function $G(z) \in \scrh(\Delta)$ a holomorphic function $H(z) \deff
\bff(\vec \zeta, \vec p, g(z)) \in \scrh(\Delta)$ of the form $H(z) = G(z) + q(z)$ such
that $q(z)$ is a polynomial of degree $\deg(q(z)) = \left(\sum_{i=1}^n m_i
\right) -1$ and such that the jets $j_{\zeta_i}^{m_i-1} h(z)$ are the given
$p_i(z)$. In other words, $h(z)$ is obtained from $g(z)$ by the prescribed
correction of its jets at the given points by means of a polynomial of the
minimal possible degree.  For example, in the case $g(z) \equiv 0$ the function $F$
realizes the Chinese remainders theorem.  It follows from the construction
that $\bff(\vec \zeta, \vec p, g(z))$ is a holomorphic function of its arguments.
The extension of $F$ to the diagonal locus of $\Delta^n$ where some of $\zeta_1, \ldots,
\zeta_n$ could coincide is maid by means of the following construction.  In a
neighborhood a point $\zeta^0=(\zeta^0_1, \ldots \zeta^0_n)$, where we have an incidence of the
form, say, $\zeta_1^0=\cdots=\zeta_k^0$, we assume that the polynomial $p_1(z), \ldots p_k(z)$
are the jets of some polynomial $\ti p(z)$ of degree $\deg(\ti p(z)) = \ti m
-1$ with $\ti m\deff \left(\sum_{i=1}^k m_i \right)$, \ie, $p_i(z) =
j_{\zeta_i}^{m_i-1} \ti p(z)$, and define $\bff(\vec \zeta^0, \vec p, g(z))$ by the
replacing the conditions on $j_{\zeta_i}^{m_i-1} h(z)$, $i= 1,\ldots,k$, by the common
condition $j_{\zeta_1^0} ^{\ti m-1} h(z) = \ti p(z)$. The function $F$ gives the
solution $q(z)$ of the equation
\begin{equation}\eqqno(1.f)
- q(z) + f(z) \prod_{j=1}^m (z - \zeta_i)^{n_i} = G(z) - \ti p(z)
\end{equation}
with unknown polynomial $q(z)$ of degree at most $\ti m -1$ and unknown $f(z)
\in \scrh(\Delta)$, in which $\zeta_i$ appear as parameters. Consequently, the regularity
of the extended $F$ is equivalent to the holomorphicity of the dependence of
$q$ (and $f(z) \in \scrh$) on the r.h.s.\ and on the parameters of the equation.
Notice also the uniqueness of such a polynomial $q(z)$.

\smallskip 
Now let $\gamma_i(t_i)$ be some irreducible holomorphic curves in $\scry_i^\circ$
connecting $R_i(C^\circ)$ with $R_i(C')$. This means that the definition domain of
each $\gamma_i$ is some irreducible curve $T_i$ and the map $\gamma_i: T_i \to \scry_i^\circ$
is holomorphic and that there exist points $t_i^\circ,t_i' \in T_i$ with $\gamma_i(t_i^\circ)
= R_i(C^\circ)$ and $\gamma_i(t_i') = R_i(C')$. By the hypotheses of the lemma we may
assume that each $\gamma_i(T_i)$ lies sufficiently close to $R_i(C^*)$. Set $T
\deff T_1 \times \cdots \times T_l$ and let $\gamma: t=(t_1, \ldots,t_l) \in T \mapsto (\gamma_1(t_1),\ldots, \gamma_l(t_l) \in
\scry_1 \times \cdots \scry_l$ be the product map. For $t=(t_1, \ldots,t_l) \in T$, let
$\{\zeta_1(t), \ldots, \zeta_n(t)\}$ be the collection of all zero points of the
discriminants $\dscr(P_{\gamma_i(t)})$ of all curves $\gamma_1(t_1), \ldots, \gamma_l(t_l)$, and
$m_1, \ldots, m_n$ the corresponding multiplicities. Then the total number $n$ of
the zeroes and their multiplicities $n_j$ are constant in $t \in T$ provided the
curves $\gamma_i(t_i)$ are chosen generic enough. Fix holomorphic map $g: T \to
\scrz_d$ with the following properties:
\begin{itemize}
\item the image $g(T)$ lies in a sufficiently small neighborhood of $C^*$;
\item the images of $t^\circ \deff (t^\circ_1, \ldots, t^\circ_l)$ and $t' \deff (t_1',
 \ldots,t_l')$ are $C^\circ$ and $C'$, respectively.
\end{itemize} 
Denote by $G_t(z,w)$ the Weierstra\3 polynomial of $g(t)$. Finally, consider
the family 
\[
H_t(z,w) \deff \bff\big( (\zeta_j(t)), 
(j^{m_j-1}_{\zeta_j(t)}P_{\gamma_{i(j)}(t_{i(j)})}),
G_t(z) \big).
\]
The meaning of the construction is as follows:
\begin{itemize}
\item We apply $\bff$ componentwisely to Weierstra\3 polynomials of the degree
 $d$, so that $H_t(z,w)$ is also a Weierstra\3 polynomial of degree $d$.
\item At each zero point $\zeta_j(t)$ of the discriminant $\dscr(P_{\gamma_i(t)})$ of
 some curve $\gamma_i(t_i)$ with the multiplicity $m_j$, we correct the $(m_j-1)$-jet
 of $G_t(z)$ to make it equal to the the $(m_j-1)$-jet
 of the Weierstra\3 polynomial $P_{\gamma_i(t_i)}$.
\end{itemize}
It follows from the construction that $H_t$ corresponds to a holomorphic map
$h: T \to \scrz_d$. Moreover, the image $h(T)$ stays close enough to $C^*$. The
uniqueness of the solution of \eqqref(1.f) implies that $H_{t^\circ} = G_{t^\circ} =
P_{C^\circ}$ and $H_{t'} = G_{t'} = P_{C'}$. Further, the condition on jets
ensures that 
\[
j^{m_j-1}_{\zeta_j(t)} \dscr(H_t) =j^{m_j-1}_{\zeta_j(t)} \dscr\big(
P_{\gamma_{i(j)}(t_{i(j)})} \big) \equiv 0.
\]
This means that the zero divisor of the discriminant $\dscr(H_t)$ is the sum
of the zero divisors of the curves $\gamma_i(t_i)$ over all $i=1, \ldots,l$. Using this
condition one can easily show that image $h(T)$ lies in $\scrz_{d, \nu}^\circ$. The
lemma follows.
\qed

\medskip
Now consider the remaining case in which \slsf{Cases 1--4} considered above
are excluded. Thus we assume that $d \geq 2$, $\nu < \delta^* \deff
\delta(C^*)$, and that the discriminant $\dscr(P_{C^*})$ vanishes only in
one point, say $z^*=0$. By \lemma{lem1.7a} we may additionally assume that
the generic curve in $\scry$ is irreducible.

We follow the idea used in \cite{Sh-2}. For each $k\in \nn$, let $F_k:
\scrz_d \to \cc^k$ be the map associating to a proper curve $C$ the
$(k-1)$-jet $j^{k-1}_0(\dscr(P_C))$ of the discriminant of its normalized
Weierstra\3 polynomial $P_C(z,w)$. Then $F_k: \scrz_d \to \cc^k$ is also
holomorphic, and the sets $\scry \cap F_k\inv(0)$ are Banach analytic of
finite definition.  Let $N$ be the order $\ord_0(\dscr(P_{C^*}))$ of the
discriminant $\dscr(P_{C^*})$ at $0$. From \eqqref(1.4) and \eqqref(1.4a) we
see that $N \geq d$ in our case. Fix a decreasing sequence $\scry= \scry_0
\supset \scry_1 \supset \cdots \scry_N$ of irreducible components $\scry_k$
of $\scry \cap F_k\inv(0)$ at $C^*$.

\newlemma{lem1.11} There exists $k^* \in \{2,\ldots d\}$ such that

\begin{enumerate}
\item[(1)] for every $k\in \{0, \ldots, k^*-1\}$ a generic curve $C \in \scry_k$ is
 has the following properties:

\begin{itemize}
\item[(1a)] $C$ is irreducible;
\item[(1b)] $C$ is nodal with $\nu$ nodes, all of them outside $\ell_0$;
\item[(1c)] $C$ has $d-k$ pairwisely disjoint non-singular branches at $\ell_0$;
\item[(1d)] the discriminant $\dscr(P_C)$ has zero of degree $k$ at $z=0$, $\nu$
 double zeroes, and $N-k -2\nu$ simple zeroes outside $z=0$.
\end{itemize}

\item[(2)] a generic curve $C \in \scry_{k^*}$ is has the following properties:
\begin{itemize}
\item[(2a)] $C$ is nodal outside $\ell_0$ with no vertical component;
\item[(2b)] the discriminant $\dscr(P_C)$ has only simple or double zeroes outside
 $z=0$; each such double zero is the projection of a node of $C$;
\item[(2c)] all local branches $C$ at $\ell_0$ are non-singular and exactly
two of them meet at $\ell_0$ whereas remaining are pairwisely disjoint from
each other and from those two; moreover, those two components either are both
vertical at $\ell_0$ or transversal to each other; 
\item[(2d)] the virtual nodal number of $C$ is $\delta(C) = \nu+1$.
\end{itemize}

\end{enumerate}
\end{lem}

\state Remark. The local behavior of a generic $C \in \scry_k$ at $\ell_0$
in the cases $(1)$ and $(2)$ is as in the cases \sli or respectively \slii
of \lemma{lem1a.3}.

\medskip
\proof Let $k_0$ be the maximal integer such that for every $k \in \{ 0,
\ldots, k_0-1\}$ a generic curve $C$ in $\scry_k$ has the properties listed
in $(1)$. Then $k_0 \geq 0$. The assertion of the lemma is that $2 \leq k_0
\leq d$ and that a generic curve $C$ in $\scry_{k_0}$ has the properties
listed in $(2)$.

We proceed using two inductive assumptions. The first one is that the
assertion of the lemma holds for any curve $C^+$ for which the order $\ord_0
\big( \dscr(P_{C^+}) \big)$ of vanishing of the discriminant
$\dscr(P_{C^+})$ at $0$ is strictly less than that for $C^*$. Another one is
a similar assumption for the multiplicity of the line $\ell_0$ in $C^+$ and
in $C^*$. The meaning of these assumptions is that the lemma holds provided
$\scry_{k_0}$ contains a curve $C^+$ which does not lie in $\scrz_d^\es
(C^*,0)$. Indeed, if the discriminant $\dscr(P_{C^+})$ of such a curve has
zeroes only at the origin $0\in \Delta$, we simply apply the lemma with
$C^*$ replaced by $C^+$. Otherwise we could apply the lemma to the curve
$C^+ \cap (\Delta(\eps) \times \pp^1)$ with $\eps>0$ small enough, making an
additional observation that a generic curve in every $\scry_k$ is nodal
outside $\ell_0$.

The remaining case when $\scry_{k_0}$ is contained in $\scrz^\es_d(C^*, 0)$.
In particular, $\codim \big( \scrz^\es_d(C^*, 0) \subset \scrz_d \big) \leq
\codim \big( \scry \subset \scrz_d \big) = \nu + k_0$ in this case. Let us
estimate $\codim \big( \scrz^\es_d(C^*, 0) \subset \scrz_d \big)$ and $\nu +
k_0$ in another way. Let $b$ be the number of irreducible non-vertical
components of $C^*$ and $m^*$ the multiplicity of $\ell_0$ in $C^*$. Note
that $C^*$ has no other vertical components. It follows from the conditions
listed in $(1)$ that $k_0 \leq d$. Using \eqqref(1.4a) we obtain the
estimate $\nu \leq \delta^* - m^* - b+1$. The meaning is that we need to
smooth at least $b+m^* -1$ nodes to obtain an irreducible curve $C \in
\scry$ from $C^*$ having $b+m^*$ components. By \lemma{lem1a.3}, $\codim
\big( \scrz^\es_d(C^*, 0) \subset \scrz_d \big) = \delta^* +m^* + d-b + e$
with $e \geq 2$ except the cases listed in the lemma. Thus
\begin{equation}\eqqno(1.e)
\delta^* + m^*+ d-b +e  \leq \nu + k_0 \leq \delta^* - m^* - b+1 +d.
\end{equation}
Consequently, $2m^* +e \leq 1$, which means $e\leq1$ and $m^*=0$. In terms
of the chosen sequence $\scry_0 \supset \scry_1 \supset \cdots \supset
\scry_{k_0}$ the condition $m^*=0$ means that the degeneration of curves $C
\in \scry$ by means of the condition $C \in \scry_{k_0}$ does not force a
``splitting out'' of a vertical component. If $e=1$ (resp.\ $e=0$), we
obtain on of the cases \slii or \sliii (resp.\ case \slip) of
\lemma{lem1a.3}. Now, we apply a new upper bound $\nu \leq \delta^* -1$.
Then we obtain
\begin{equation}\eqqno(1.g)
\delta^* + d-b +e  \leq \nu + k_0 \leq \delta^* -1 +d,
\end{equation}
and, consequently, $b \geq e+1$.  This excludes case \sliii of \lemma{lem1a.3}
since in our situation curve $C^*$ must be connected which would imply $b=1$
in contradiction with $b\geq e+1 =2$. Observe also that $e=1$ means that we must
have the equalities in \eqqref(1.e), and in particular, $\nu = \delta^* -1$. In the
considered special situation, when $\scry_{k_0} \subset \scrz^\es_d(C^*, 0)$, the
latter is equivalent to the condition (2d).
\qed

\medskip \proofp of \refthm{thm1.9}. Let $C^* \in \scrz_d$ be a proper curve and
$\scry$ a component of $\scrz_{d,\nu}$ passing through $C^*$. Applying the
previous lemma, we reduce the general situation to the case when $C^*$ has the
properties (2a--2d) of the lemma. Let $p$ be the point on $\ell_0$ such that
there are two local branches of $C^*$ at $p$. Then there exist a
neighborhood $U$ of $p$ and complex coordinates $(\ti z, \ti w)$ in $U$ with
the following properties:
\begin{itemize}
\item there exists a biholomorphic map $\phi:U \xarr{\;\cong\;}{} \Delta^2$ which
 extends holomorphically into some neighborhood of the closure $\barr U$;
\item $(\ti z, \ti w)$ are the pull-back of the the standard coordinates in
 $\Delta^2$ with respect to $\phi$;
\item there exists a neighborhood $\scru \subset \scrz_d$ of $C^*$ such that for
 every curve $C \in \scru$ the curve $C \cap U$ is defined by a Weierstra\3
 polynomial $\ti w^2 + \ti a_1(\ti z) \ti w + \ti a_2(\ti z)$ of degree 2;
\item the discriminant of the Weierstra\3 polynomial of the curve $R_U(C^*)$
 vanishes only at the origin $\ti z =0$.
\end{itemize}
It can be easily seen that the map $R_U : \scru \to \scrz_2(U)$ given by $C \in
\scru \mapsto C \cap U \in \scrz_2(U)$ is holomorphic, and that each preimage
$R_U\inv(\scrz_{2, \nu}(U))$ is a union of some components of $\scrz_{d,\nu}(\Delta \times
\pp^1)$.  Set $\delta_p \deff \delta(C^*, p)$ and $\nu_p \deff \delta_p-1$. Then for a generic
curve $C \in \scry$ close enough to $C^*$ the curve $R_U(C)$ is nodal with $\nu_p$
nodes. The crucial point in the proof of \refthm{thm1.9} is the following
assertion:

\begin{itemize}
\item[\ ] \sl Every local irreducible component $\scrw$ of $\scrz_{2, \delta_p}(U)$ at
 $R_U(C^*)$ contains $\scrz_{2, \nu_p}(U)$.
\end{itemize}

The assertion is the special case of \refthm{thm1.9} for the curve $R_U(C^*)$.
To prove it, we simply apply \lemma{lem1.11} to the component $\scrw$ and
obtain the locus $ \scrw_2 \subset \scrw$ such that a generic curve in $\scrw_2$ is
nodal with $\nu_p +1 = \delta_p$ nodes. 

To deduce the proof of the theorem, let us consider the locus $\scry' \deff
\scry \cap R_U\inv(\scrz_{2,\delta_p})$. Then $\scry'$ is a BASFD, $C^*$ lies in
$\scry'$. Further, a generic curve $C \in \scry'$ close enough to $C^*$ is nodal
with $\nu+1$ node. The additional node appears in the neighborhood $U$ of the
point $p$. Repeating this construction we can produce one by one all possible
additional nodes on a curve in $\scry$ in such a way that the curve will
remain nodal. The procedure stops when we achieve the locus of maximal nodal
deformations of $C^*$.
\qed

\bigskip

\newsection[glob]{Severi problem for Hirzebruch surfaces} 

\newsubsection[glob1]{Severi problem for ruled surfaces} We recall briefly the
definition and main properties of ruled surfaces. 

\newdefi{def2.1} A \slsf{ruling} of a smooth complex surface $X$ over a smooth
complex curve $Y$ is a proper holomorphic projection $\pr: X \to Y$, such that
$d\pr : T_xX \to T_{\pr(x)}Y$ is surjective for every $x\in X$ and such that for
generic $y\in Y$ the fiber $\pr\inv(y)$ is isomorphic to the complex projective
line $\pp^1$. In this case $X$ is called a \slsf{ruled surface} and $Y$ the
\slsf{base of the ruling}. We consider the ruling $\pr: X \to Y$ as a part of
the structure of $X$. A fiber $\pr\inv(y)$ of a ruling $\pr: X \to Y$
isomorphic to $\pp^1$ is called \slsf{regular} or a \slsf{vertical line} and
denoted by $\ell_y \deff \pr\inv(y)$; a non-regular fiber is called
\slsf{singular}.

A ruling $\pr: X \to Y$ is called \slsf{minimal} if every fiber is regular.
In this case $X$ is called a \slsf{minimal ruled surface}.

In most case we assume that such a surface $X$ (and hence the base $Y$) is
compact.
\end{defi}

The structure of ruled surfaces is well understood, so we only list some its
properties referring to the standard sources \cite{B-P-V, Gr-Ha, Hart} for a
more detailed exposition.

A (non-singular) compact complex surface $X$ admits a ruling $\pr: X \to Y$ iff
there exists a non-singular rational curve $C \subset X$ with self-intersection $C \cdot
C = 1$. In this case $c_1(X)\cdot C =2$ by genus formula and the ruling of $X$ can
be constructed as the family of deformations of $C$ on $X$; in particular,$C$
is a fiber of this ruling.  Every non-minimal ruled surface $X$ can be
obtained as a blow-up of a minimal ruled surface $X'$ such that the ruling
$\pr_X: X \to Y$ is the composition of the contraction map $\pi: X \to X'$ with the
ruling $\pr_{X'} : X' \to Y$. Such a contraction $\pi: X \to X'$ is always not
unique.  For example, blowing-up a point on a regular fiber of a ruled surface
$X$ with the projection $\pr: X \to Y$ we obtain a singular fiber consisting of
two exceptional curves; the first one, say $C'$, is the exceptional curve of
the blow-up and the other, say $C''$, is the proper pre-image of the fiber
containing the center of the blow-up. Contracting $C''$ we obtain a new 
non-singular complex surface on which the original exceptional curve $C'$
becomes a regular fiber. 

A non-compact minimal ruled surface is isomorphic to $Y\times\pp^1$ and the
projection $\pr: Y \times\pp^1$ is its unique possible ruling. Every compact minimal ruled
surface except the blown-up $\pp^2$ is minimal as an abstract complex surface.
Every compact minimal ruled surface except $\pp^1 \times \pp^1$ has a unique
ruling. Every compact minimal ruled surface $X$ with the ruling $\pr: X \to Y$
has the form $X = \pp(E)$ where $E$ is some holomorphic vector bundle over $Y$
of rank $2$, and the ruling $\pr: X \to Y$ is induced by the projection $E \to Y$.
In particular, every compact ruled surface is projective. Two holomorphic
vector bundles $E_1, E_2 \to Y$ define isomorphic surfaces $\pp(E_1),\;
\pp(E_2)$ iff $E_1 \cong E_2 \otimes L$ for some holomorphic line bundle $L \to Y$. The
surfaces $\pp(E)$ and $\pp(E^*)$ --- where $E^* \deff \scrhom(E, \scro_Y)$ is
the dual bundle --- are isomorphic.  This fact follows immediately from the
isomorphism $E^* \cong E \otimes \big( \det(E) \big)\inv$.

For a given compact minimal ruled surface $X$ with the ruling $\pr: X \to Y$
there exists a rank 2 vector bundle $E$ over $Y$ such that $\sfh^0(Y, E) \neq 0$
but $\sfh^0(Y, E \otimes L) = 0$ for any holomorphic line bundle $L$ on $Y$ of
negative degree $c_1(L)$. We call such a bundle $E$ a \slsf{normalized vector
 bundle} defining $X$. For a given $X$ there exist finitely many normalized
vector bundles $E$ defining $X$ and all of them have the same degree $c_1(E)$.
The number $e \deff - c_1(E)$ is called \slsf{Hirzebruch index $e(X)$} of the
minimal ruled surface $X$. A normalized vector bundle $E$ is the extension of
the form $0 \to \scro_Y \to E \to Q \to 0$ with $Q = \det(E)$.  A compact minimal
ruled surface $X$ of the $X = \pp(E)$ is of \slsf{split type} if so is its
defining vector bundle $E$, \ie $E \cong \scro_Y \oplus \det(E)$. This is equivalent to
the splitting $E' \cong L_1 \oplus L_2$ of any holomorphic vector bundle $E'$ defining
$X$, \ie, each time when $\pp(E')= X$. The Hirzebruch index $e(X)$ is always
non-negative for ruled surfaces of split type and varies in the range $g_Y \leq
e(X) \leq 2g_Y-2$ for ruled surfaces $X$ over a curve $Y$ of genus $g_Y$.
Moreover, any $e \geq0$ (resp., any $e$ in the range $g_Y \leq e \leq 2g_Y-2$) is
realizable by an appropriate minimal ruled surface of (non-)split type.

The Severi problem for ruled surfaces can be formulated as follows:
\begin{quote}
 \it Describe connected components of the locus $\scrz_\nu^\circ(X, [A])$ of nodal
 curves in on a given compact ruled surface $X$ with a given nodal number $\nu$
 and given homology class $[A] \in \sfh_2(X, \zz)$. What are necessary or
 sufficient conditions for irreducibility of $\scrz_\nu(X, [A])$?
\end{quote}
Below we give some examples of situations where $\scrz_\nu^\circ(X, [A])$ has
expected dimension 0 and consists of several points. This leads to a more
sophisticated version of the Severi problem asking whether one can pass from
one irreducible component of $\scrz_\nu^\circ(X, [A])$ to another one using the
monodromy of some family $\{X_s\}_{s \in S}$ of deformations of the given surface
$X$. Therefore it is interesting to determine what complex structures on
compact ruled surfaces are ``generic''. To give a precise sense to this
notion, let us recall that the deformation theory (see \eg \cite{Pal-1,
 Pal-2}) provides a semi-universal family $\{\scrx_s\}_{s\in S}$ of deformations
of a given compact complex manifold $X$ whose base $S$ can be realized as an
analytic set in a ball in the space $\sfh^1(X, \scro^{TX})$. As above, we say
that some property $\frA$ is (Zariski-analytic) \slsf{generic} for a given
class of compact complex manifolds if for any manifold $X$ in this class there
exists an analytic set $S_\frA$ of the base $S$ of semi-universal family
$\{\scrx_s\}_{s\in S}$ of deformations of $X$ such that $A$ does not contain any
irreducible component of $S$ and such that the some property $\frA$ holds for
any $\scrx_s$ with $s \in S \bs S_\frA$. Moreover, as such an analytic set
$S_\frA$ one can take the Zariski-analytic closure of the locus $S_\frA^\circ$ of
$s \in S$ parameterizing those deformations $\scrx_s$ of $X$ which have the
property $\frA$. Here we assume implicitly that the class of complex manifolds
we consider is stable under deformations.

\newlemma{lem2.0}
\begin{itemize}
\item[\sli] A {\sl generic} minimal compact ruled surface $X$ with the ruling
 $\pr: X \to Y$ over a curve $Y$ of genus $g$ has Hirzebruch index $e(X) = -(g-1)$
 or $e(X) = -g$. Moreover, $X$ is of split type if $e(X) \geq 0$ and of non-split
 type otherwise.
\item[\slii] Every singular fiber of a {\sl generic} non-minimal compact ruled
 surface $X$ is a union of two exceptional rational curves meeting
 transversally at a single point. 
\end{itemize}
\end{lem}

\newdefi{def2.2} A fiber of a ruling $\pr:X \to Y$ consisting of two exceptional
rational curves meeting transversally at a single point is called an
\slsf{ordinary singular fiber}.
\end{defi}

Before giving the proof, let us now describe curves on a given compact minimal
ruled surface $X$ with a fixed normalized defining vector bundle $E$ over a
curve $Y$. Denote by $\pr_E: E \to Y$ the projection map. There exists an open
covering $U_\alpha$ of $Y$ such that each $E_\alpha \deff \pr_E\inv(U_\alpha)$ is
(holomorphically) isomorphic to $U_\alpha \times \cc^2$. This gives us local coordinates
$z_\alpha \in U_\alpha$ and $w_\alpha = [w_{\alpha,0}: w_{\alpha,1}] \in \pp^1$ on each $X_\alpha \deff
\pr\inv(U_\alpha) \cong U_\alpha \times \pp^1$.  For a curve $C$ in $X$ we define its degree $d$
over $Y$ as the intersection index $C \cdot \ell_y$ with any vertical line $\ell_y =
\pr\inv(y)$. In every chart $X_\alpha$ the curve $C$ in defined by a Weierstra\3
polynomial $P_{C,\alpha}(z_\alpha, w_\alpha) = \sum_{i=0}^d a_{i, \alpha}(z_\alpha)w_{\alpha,0}^{d-i}
w_{\alpha,1}^i$. We consider each $P_{C,\alpha}$ as a section of the symmetric power
$\sym^d(E^*)$ over $U_\alpha$. Then $P_{C,\alpha} = g_{\alpha\beta} \cdot P_{C,\beta}$ for some
holomorphic non-vanishing functions $g_{\alpha\beta} \in \scro^*(U_{\alpha\beta})$ where $U_{\alpha\beta}
\deff U_\alpha \cap U_\beta$. So the system $\{g_{\alpha\beta}\}$ form a cocycle, \ie $g_{\alpha\beta} \cdot
g_{\beta\gamma} = g_{\alpha\gamma}$ in $U_{\alpha\beta\gamma} \deff U_\alpha \cap U_\beta \cap U_\gamma$, such that $P_{C,\alpha}$ can
be considered as local trivializations of the line bundle $L$ given by the
cocycle $\{g_{\alpha\beta}\}$. On the other hand, since every $P_{C,\alpha}$ is a section of
$\sym^d(E^*)$, we obtain a holomorphic homomorphism of bundles $F: L \to
\sym^d(E^*)$. Vice versa, for any holomorphic line bundle $L$ and any non-zero
homomorphism $F: \to \sym^d(E^*)$ we consider a system of local trivializations
$s_\alpha \in \scro(U_\alpha)$ of $L$ over an open covering $\{U_\alpha\}$ of $Y$ and set $P_\alpha
\deff F(s_\alpha) \in \sfh^0(U_\alpha, \sym^d(E^*))$. Then $P_\alpha$ form a defining system of
Weierstra\3 polynomials of a curve $C$ of degree $D$ on $X$. In this way we
come to

\newlemma{lem2.2} There exists a 1-to-1 correspondence between curves $C$ of
degree $d$ on a compact minimal ruled surface $X$ with a fixed defining vector
bundle $E$ over a curve $Y$ and coherent subsheaves $\scrs$ of
$\rank(\scrs)=1$ in symmetric powers $\sym^d(E)$. Moreover, if a rank $1$
subsheaf $\scrs \subset \sym^d(E)$ corresponds to a curve $C \subset X$, then its
saturation $\scrs^{\perp\perp} \subset \sym^d(E)$ corresponds to the union of all
non-vertical components of $C$. Furthermore, the curve corresponding to a
saturated rank $1$ subsheaf $\scrs \subset E$ is the projectivization $\pp(\scrs)
\subset \pp(E) = X$ of $\scrs$.
\end{lem} 

Here, in abuse of notation, we identify holomorphic bundle $E$ with
the sheaf $\scro(E)$ of its holomorphic sections denoting the latter also by
$E$.

The lemma allows to give a pure algebraic definition of the spaces $\scrz_\nu^\circ
(X, [A])$ for compact ruled surfaces. Namely, in the case of minimal $X =
\pp(E)$ we the spaces $\scrz^\circ_\nu(X, [A])$ can be described as algebraic sets
in the locus $\scrz(X, [A])$ of the pairs $(L, [s])$ where $L$ is a line
bundle over the base curve $Y$ of a certain degree $a$ and $[s]$ a point of
$\pp(\sfh^0(Y, \sym_dE \otimes L))$ with an appropriate $d$. In the case of
non-minimal ruled surface $X$ the description involves the images of the
curves under the projection  $f: X \to X'$ onto an appropriate minimal ruled
surface $X'$.

\proof Recall that the saturation of a subsheaf $\scrs$ of a torsion free
sheaf $\scre$ is the subsheaf $\ti \scrs \subset \scre$ which contains $\scrs$ and
such that the quotient $\scre/ \ti\scrs$ is torsion-free. Furthermore, if
$\scre$ is a reflexive sheaf, \eg, a locally free sheaf, then the saturation
can be constructed as the double orthogonal $\scrs^{\perp\perp}$ of $\scrs$. Every
reflexive sheaf $\scre$ on a curve is locally free, \ie, is the sheaf of local
(holomorphic) sections of the uniquely defined vector bundle $E$, and the
saturation $\scrs^{\perp\perp}$ of a subsheaf $\scrs \subset \scre$ on a curve is a
subbundle of $E$. On a curve $Y$, a rank 1 subsheaf $\scrs$ of the locally
free sheaf $E^*$ is locally generated by one section. Such a section $s$
admits a local representation $s(z) = h(z) \cdot s'(z)$ with a non-vanishing
section $s'$ of $\scre$ and a local holomorphic function $h(z)$. Considering
$s$ as a local equation of a curve $C\subset X$, we see that the decomposition $s=
h\cdot s'$ induces the local decomposition $C= C' \cup \bigcup_i{\ell_{z_i}}$ where $C'$ is a
curve with no vertical components and $\{z_i\}$ is the divisor of $h$. To finish
the lemma, it remains to use the natural isomorphism $E^* \cong E \otimes \det(E^*)$
which is valid for all (holomorphic) vector bundles of rank 2.
\qed

\proofp of \lemma{lem2.0}. We make use of the following sufficient condition
for genericity of a given property $\frA$. If for any compact complex surface
$X$ there exists a family $\{\scrx_s\}_{s\in S}$ of deformations of $X$ with an
irreducible base $S$ and a proper analytic subset $S_\frA \subset S$ such that
$\frA$ holds for any $\scrx_s$ with $s \in S \bs S_\frA$, then $\frA$ holds for
a generic compact complex surface (in a given class). Indeed, in this case the
the locus of deformations of $X$ without property $\frA$ can not be
Zariski-analytic dense in any component of the semi-universal family of
deformations of $X$.

\slsf{Part \sliip} follows in view of this sufficient condition rather easily.
Indeed, since each non-minimal ruled surface $X$ is a blow-up of some minimal
ruled surface $X'$. The corresponding blow-up center is, in general, a zero
dimensional non-reduced subspace $Z$ of $X'$ of a given length of the
structure sheaf $\scro_Z$.  Due to Douady \cite{Dou}, there exists a
holomorphic family of deformations $\{Z_s\}$ of such $Z$ parameterized by a
complex analytic space $S$.  Moreover, since $X'$ is algebraic, such a
parameterizing space is an analytic chart of the appropriate Hilbert scheme of
points on $X'$. Blowing-up $X'$ in $Z_s$ we obtain a holomorphic family
$\{\scrx_s\}$ of deformations of $X$ parameterized by the same space $S$. An
easy observation is that $\scrx_s$ has only ordinary singular fibers iff each
fiber of $X'$ is blown-up at most once. This means that $Z_s$ has no multiple
points, \ie, the length of each local ring $\scro_{Z_s, p}$ is $1$, and that
each fiber $\ell'_y$ of the ruling of $X'$ contains at most $1$ point from $Z_s$.
So it remains to notice that such a situation holds for a generic $s \in S$.
Finally, observe that the argument works as well in the case of non-compact
ruled surfaces.

\medskip
\slsf{Part \slip.} Let $X$ be a minimal compact ruled surface of the form $X = \pp(E)$
with a holomorphic vector bundle $E$ over a curve $Y$. The remark above and
\lemma{lem2.2} allow us to reduce the first assertion of the lemma to the
problem of the (non-)existence of global holomorphic section of the bundles
$E' \otimes L$ where $E'$ is some deformation of $E$ and $L$ a holomorphic line
bundle of given degree. 

By the theorem A, for an appropriate line bundle $L_1$ of a sufficiently high
degree the bundle $E \otimes L_1$ admits a non-vanishing section $s \in \sfh^0(Y, E \otimes
L_1)$. This gives us the extension 
\begin{equation}\eqqno(2.1)
0 \to \scro \xarr{s} E \otimes L_1 \to Q \to 0
\end{equation} 
where $\scro$ is the trivial line bundle over $Y$ and $Q$ the quotient line
bundle. Thus we can include $E \otimes L_1$ into the family $\{E_\xi\}$ where $E_\xi$ is
the extension 
\begin{equation}\eqqno(2.2)
0 \to \scro \to E_\xi \to Q \to 0
\end{equation}
with a given $\xi \in \sfh^1(Y, Q^*) \cong \sf{Ext}^1(Y; Q, \scro)$. Projectivizing,
we obtain a holomorphic family $\scrx_\xi \deff \pp(E_\xi)$ of deformations of $X$.

Notice that for any line bundle $L \in\sf{Pic}(Y)$ the ``twisted'' sequence $0
\to L\to E_\xi \otimes L\to Q \otimes L\to 0$ is also induced by $\xi$ by means of the isomorphism
$\sfh^1(Y, Q^*) \cong \sf{Ext}^1(Y; Q, \scro) \cong \sf{Ext}^1 (Y; Q \otimes L, \scro \otimes L)$.
Moreover, the corresponding connecting homomorphism $\delta: \sfh^0(Y, Q \otimes L) \to
\sfh^1(Y, L)$ is also given as the product with $\xi$ with respect to the Yoneda
multiplication
\[
\sfh^0(Y, Q \otimes L) \otimes \sfh^1(Y, Q^*) \mapsto \sfh^1(Y, Q \otimes L \otimes Q^*) = \sfh^1(Y, L).
\]
We use the notation $\xi_*$ to denote the maps induced by the Yoneda
multiplication with $\xi_*$.

Denote by $l$ the degree of $L^*$ and by $q$ the degree of $Q$. By the
construction, $Q = \det(E \otimes L_1)$ has sufficiently high degree, $q \gg0$. In the
situation we are interested in $L$ is a negative line bundle, which means that
$l >0$. In this case every section of $E_\xi \otimes L$ descents to a section $\gamma$ of
$Q \otimes L$ such that $\xi_*(\gamma)=0 \in \sfh^1(Y,L)$, and vice versa. As it is done with
$\xi$, we denote by $\gamma_*: \scro \to Q \otimes L$ the bundle homomorphism induced by $\gamma \in
\sfh^0(Y,Q \otimes L)$. In this way we come to the commutative diagram
\begin{equation}\eqqno(2.3)
\xymatrix{ 
\sfh^0(Y,\scro) \ar[r]^{\gamma_*} \ar[d]^{\xi_*}& \sfh^0(Y, Q \otimes L) \ar[d]^{\xi_*}\\
\sfh^1(Y,Q^*) \ar[r]^{\gamma_*} & \sfh^1(Y, L).
}
\end{equation}
Dualizing it, we obtain 
\begin{equation}\eqqno(2.4)
\xymatrix{ 
\sfh^1(Y,K) \ar@{<-}[r]^{\gamma_*} \ar@{<-}[d]^{\xi_*}& \sfh^1(Y, K \otimes Q^* \otimes L^*) \ar@{<-}[d]^{\xi_*}\\
\sfh^0(Y,K \otimes Q) \ar@{<-}[r]^{\gamma_*} & \sfh^0(Y, K \otimes L^*),
}
\end{equation}
where $K$ denoted the canonical line bundle on $Y$, $K\deff \Omega^1_Y$. Since $l =
\deg(L^*) >0$ by our assumption, the space $\sfh^0(Y, K \otimes L^*)$ has dimension
$g-1 +l$ where $g$ is the genus of $Y$. Similarly, $\dim \sfh^0(Y, K \otimes Q) =
g-1 +q$. Assume that $\xi\neq 0$ and denote by $W_\xi$ the kernel of the map $\xi_*:
\sfh^0(Y, K \otimes Q) \to \sfh^1(Y,K) \cong \cc$ and by $[\xi]$ the corresponding point in
$\pp((\sfh^0(Y, K \otimes Q))^*) \ddef \pp^{g-2+q}$. Finally, denote by $D_\gamma$ the
divisor of $\gamma$ and by $D_\sigma$ the divisor of a given $\sigma\in \sfh^0(Y,K \otimes Q)$.
Consider $S_\gamma \deff \gamma_*(\pp(\sfh^0(Y, K \otimes L^*)))$ as a linear subsystem in
$\pp(\sfh^0(Y, K \otimes Q))$. Then the existence of $\gamma\in \sfh^0(Y, Q \otimes L)$ with
$\xi_*(\gamma)=0$ is equivalent to the existence of an effective divisor $D_\gamma$ of
degree $q-l$ and a linear subsystem $S$ in the linear system $\pp(W_\xi) \subset
\pp(\sfh^0(Y, K \otimes Q))$ of dimension $\dim(\pp(\sfh^0(Y, K \otimes L^*))) = g-2+l$
such that $D_\sigma \geq D_\gamma$ for any $\sigma\in S$. In terms of the corresponding imbedding
$\phi_{K\otimes Q}: Y \to \pp((\sfh^0(Y, K \otimes Q))^*) = \pp^{g-2+q}$ we obtain the
following interpretation: There exists a linear space $S^\perp \subset \pp^{g-2+q}$ of
dimension $(g-2+q) - (g-2+l)-1=q-l-1$ which passes through the point $[\xi]$
(condition $S \subset W_\xi$) and through all points of $\phi_{K\otimes Q}(D_\gamma)$ (condition
$D_\sigma \geq D_\gamma$). The latter condition is interpreted in the usual sense in the
case when $D_\gamma$ has multiple points. Namely, if $D_\gamma = \sum m_iy_i$ with $m_i\geq1$
and $y_i \in Y$, then $\phi_{K\otimes Q}(y_i) \in S^\perp$ and $\phi_{K\otimes Q}(Y)$ has osculation
with $S^\perp$ of order $m_i -1$ at $\phi_{K\otimes Q}(y_i)$. Since $\gamma$ is a section of the
line bundle $Q \otimes L$, the degree of its divisor $D_\gamma$ is $q-l$. Thus for a
generic choice of the divisor $D_\gamma$ $S^\perp$ must be a $(q-l-1)$-plane spanned by
$\phi_{K\otimes Q}(D_\gamma)$. The variety of points swept by all such $(q-l-1)$-plane
corresponding to all possible divisors $D$ of degree $q-l$ has dimension $q-l
+ q-l-1 = 2q-2l-1$.  Thus in the case $g-2+q > 2q-2l-1$ a generic $[\xi] \in
\pp^{g-2+q}$ is not contained in any $(q-l-1)$-plane $S^\perp$ which can be
spanned by $q-l$ points lying on $\phi_{K\otimes Q}(Y)$. Taking $\xi \in \sfh^1(Y, Q) \cong
(\sfh^0(Y, K \otimes Q))^*$ with this property and the corresponding extension
$E_\xi$, we obtain a ruled surface $X_\xi = \pp(E_\xi)$ such that the bundle $E_\xi \otimes
L$ has no section for any line bundle $L$ of degree $-l$ satisfying $g-2+q >
2q-2l-1$. Since $q = \deg(E_\xi)$ and $g-2+q > 2q-2l-1$ is equivalent to $g-1 >
q-2l$, the degree of the normalized vector bundle representing $X_\xi$ is at
least $g$, \ie, $e(X_\xi) \leq -(g-1)$. On the other hand, as it was already
noticed $e(X) \geq -g$. So $e(X_\xi) = -g$ and $e(X_\xi) = -(g-1)$ are the only
possibilities.

As it was already noticed, the minimal ruled surface $X$ has non-split type in
the case $e(X) <0$. In this case the genus $g$ of the base of $X$ is $0$ or
$1$. Since every holomorphic vector bundle on $\pp^1$ splits by the classical
theorem of Grothendieck, it remains to consider the case $g=1$ and $e(X)=0$.
Let $X$ be a minimal ruled surface $X$ of non-split type over a a base curve
$Y$ of genus $g=1$ such that $e(X)=0$. Then $X$ has the form $X = \pp(E)$ for
some vector bundle $E$ which can be included in a non-trivial extension $0\to
\scro \to E \to L \to 0$ with some line bundle $Q$ of degree $\deg(L)=0$. If $L \not
\cong \scro$, then $\sfh^1(Y, L^*) \cong (\sfh^0(Y, L))^* =0$, so that $E$ must split
in this case, in contradiction to our assumption. Thus $L \cong \scro$ in our case.
For a divisor $D= \sum_i m_iy_i$ on $Y$ and a holomorphic vector bundle $F$ over
$Y$, denote by $\scro[D]$ the line bundle associated with the divisor $D$ and
by $F[D]$ the bundle $F \otimes \scro[D]$. Considering the induced sequence $0\to
\scro[y] \to E[y] \to \scro[y] \to 0$ it is easy to see that $\sfh^0(Y, E[y])$ has a
basis $\{s_0, s_1\}$ such that $s_0$ generates $\sfh^0(Y, \scro[y]) \subset \sfh^0(Y,
E[y])$ and vanishes at $y$. We observe that $s_1$ never vanish. Indeed, if
$s_1(y)=0$, then $\dim\sfh^0(Y, E) \geq 2$ which would imply the splitting of
$E$. Similarly, if $s_1(y')=0$ with some $y'\neq y$, then $\dim \sfh^0(Y,
E[y-y']) \geq1$ which would contradict $\sfh^0(Y, \scro[y-y']) =0$. Thus the
quotient $E[y] / s_1\scro$ is a line bundle which is isomorphic to $\det(E[y])
=\scro[2y]$. Consider the obtained extension $0 \to \scro \xarr{s_1} E[y] \to
\scro[2y] \to 0$ and the associated connecting homomorphism $\delta: \sfh^0(Y,
\scro[2y]) \to \sfh^1(Y,\scro) \cong \cc$. Then the section $s_0$ of $E[y]$ projects
to a non-trivial section $s$ of $\scro[2y]$ such that $s(y)=0$ and $\delta(s)=0$.
This means that $s$ is generates the image of the natural imbedding $\sfh^0(Y,
\scr0) \to \sfh^0(Y, \scro[2y])$. Summing up we obtain the following
characterization of minimal ruled surfaces $X$ over an elliptic curve $Y$ which
have non-split type and Hirzebruch index $e(X)=0$. Every such $X$ has the
form $X=\pp(E)$ with some rank $2$ holomorphic vector bundle $E$ on $Y$ which
can be included into the extension $0 \to \scro \to E \to \scro[2y] \to 0$ with any
given $y \in Y$ such that the kernel of the corresponding connecting
homomorphism $\delta: \sfh^0(Y, \scro[2y]) \to \sfh^1(Y,\scro) \cong \cc$ is the space
$\sfh^0(Y, \scro) \to \sfh^0(Y, \scro[2y])$. On the other hand, since the
canonical bundle of any elliptic curve is trivial, we can identify the
homomorphism $\delta$ with the element $\xi \in \sf{Ext}^1(Y; \scro[2y], \scro) \cong
\sfh^1(Y, (\scro[2y])^* ) \cong \hom(\sfh^0(Y, \scr[2y]), \cc)$. This implies that
for a generic $\xi \in \sfh^1(Y, (\scro[2y])^* )$ the extension $0 \to \scro \to E_\xi \to
\scro[2y] \to 0$ defined by $\xi$ yields a vector bundle $E_\xi$ of split type,
hence such is the surface $X_\xi \deff \pp(E_\xi)$. Thus we obtain a family
$\{X_\xi\}_{\xi\in \sfh^1(Y, (\scro[2y])^* )}$ of deformations of $X$ whose generic
member is of split type. The lemma follows.
\qed

\newsubsection[sec2.2]{Sections of ruled surfaces} We use
the technique developed in the proof of \lemma{lem2.0} to count the number of
the curves $C$ of degree $1$ with $[C]^2 = g-1$ on a generic compact minimal
ruled surface over a base curve $Y$ of genus $g$. By \lemma{lem2.2}, each
irreducible curve $C$ of degree $1$ on a ruled surface $X$ of the form $X=
\pp(E)$ corresponds to a line subbundle $L \subset E$. The ruling projection $\pr: X
\to Y$ maps $C$ isomorphically onto $Y$. Thus $C$ defines a section $\sigma: Y \to X$
of the projection $\pr: X \to Y$ such that $C= \sigma(Y)$. Let us consider the
extension $0 \to \scro \to E \otimes L\inv \to L_1 \to 0$. Lifting this extension to $C$ and
using the equality $\pp(E \otimes L\inv) =X$, we see that $L_1$ (more precisely, to
its lift $\sigma_*(L_1)$ from $Y$ to $C$) is isomorphic no the normal bundle to
$C$. Another proof of this fact can be obtained from the adjunction formula
combined with the formula for the canonical class of ruled subcases, see \eg
\cite{Hart}, {\S}\,V.2. In particular, the degree of the normal bundle $N_C$ is
$\deg(E \otimes L\inv)= \deg(E) -2 \deg(L)$. Since $C$ is imbedded, $\deg(N_C) =
[C]^2 = -(g-1)$. Assume that $E$ is normalized, so that it can be included
into the extension $0 \to \scro \to E \to Q \to 0$ with $\deg(Q) = -e(X)$. Since $X$
is generic, $\deg(Q) = g-1$ or $\deg(Q) = g$. In the latter case $[C]^2 \not \equiv
g-1 \mod2$ for any curve of degree $1$. So we consider the non-trivial case
$\deg(Q) = g-1$.

Thus we are interested in the number of line subbundles $L$ in a given generic
bundle $E$ over a given curve $Y$ of genus $g$ which can be included into the
extension 
\begin{equation}\eqqno(2.5)
0 \to \scro \to E \to Q \to 0
\end{equation}
with $\deg(Q) = g-1$. One such subbundle is $\scro \subset E$. In the case $g=0$ we
have $Y \cong \pp^1$. Thus the bundle $E$ splits into the sum $E \cong \scro \oplus
\scro(-1)$, and $\scro$ is the unique subbundle with the desired properties.
In the case $g=1$ the bundle $E$ also splits into the sum $E \cong \scro \oplus Q$.
Moreover, $Q \not \cong \scro$ for generic $E$. It is easy to show that $\scro$
and $Q$ are the only line subbundles of $E$ of degree $0$.  Now consider the
family $X_Q \deff \pp (\scro \oplus Q)$ of deformations of $X$ where $Q$ varies in
the Picard variety $\sf{Pic}_0(Y) \bs \{0\}$ of non-trivial line bundle on $Y$
of degree $0$. Then on each $X_Q$ we obtain two curves $C_0$ and $C_1$ of
degree $1$ with $[C_i]^2 =0$ corresponding to the summands $\scro$ and $Q$,
respectively. Denote by $Q_0$ the parameter corresponding to the original
surface $X$ and set $Q_1 \deff Q \inv$. Then $X_{Q_1}$ is isomorphic to
$X_{Q_0} = X$, since the bundles $\scro \oplus Q_0$ and $\scro \oplus Q_1$ differ by
multiplication with a line bundle, $\scro \oplus Q_0 \cong (\scro \oplus Q_1) \otimes Q_0$. But
the latter isomorphism interchanges the summands. Consequently, the isomorphism
between $X_{Q_0}$ and $X_{Q_1}$ interchanges the curves $C_0$ and $C_1$. Thus
the monodromy along the family $\{X_Q\}$ acts transitively on curves $C_0,\;
C_1$ on $X$. 

\smallskip 
In the case $g=2$ the surface $X$ is represented by a non-split vector bundle
$E$ admitting the extension \eqqref(2.5) with $\deg(Q)=1$. Let $\xi$ be the
element of $\sf{Ext}^1(Y; Q, \scro) \cong \sfh^1(Y, Q\inv) \cong \big( \sfh^0(Y, K \otimes
Q) \big)^*$ defining the extension \eqqref(2.5) and $L \subset E$ a subbundle of $E$
of degree $0$ different from $\scro \subset E$. It was shown in the proof of
\lemma{lem2.0} that to every line subbundle $L \subset E$ of degree $0$ different
from $L_0 \deff \scro \subset E$ we can associate an effective divisor $D_\gamma$ of
degree $g-1$ such that there exists a divisor $D_\sigma \geq D_\gamma$ in the linear system
given by the space $W_\xi = \ker( \xi: \sfh^0(Y, K \otimes Q) \to \cc)$. In our case $g=2$
the space $W_\xi$ has dimension $g-1=1$, so the linear system $\pp(W_\xi)$
consists of a unique divisor $D_\sigma$. Its degree is $\deg(K \otimes Q) = 3g-3 =3$. The
degree of $D_\gamma$ is $g-1=1$, so $D_\gamma$ is one of the points of $D_\sigma$. For a
generic choice of $\xi$ the divisor $D_\sigma$ consists of 3 pairwise distinct point.
The possibility to invert the construction insures that each of 3 points of
$D_\sigma$ yields a line subbundle $L \subset E$ with the desired properties.

To study the action of the monodromy group on the curves $C_0,C_1, \ldots, C_3$
corresponding to the constructed line bundles $L_0\cong \scro, L_1, \ldots ,L_3$ we
consider the locus $S$ of triples $\sigma=\{y_1, y_2, y_3\}$ of points of the base
curve $Y$ with pairwise distinct $y_1 \neq y_2 \neq y_3 \neq y_1$.  For every such $\sigma \in
S$ we set $D_\sigma \deff y_1+y_2+y_3$, $Q_\sigma \deff \scro[D_\sigma] \otimes K\inv$ and define
$W_\sigma$ to be the space of sections of $\scro[D_\sigma] = K \otimes Q_\sigma$ with zero divisor
$D_\sigma$. Then $\dim\sfh^0(Y, K \otimes Q_\sigma)= 2$ and $W_\sigma$ is a subspace of $\sfh^0(Y,
K \otimes Q_\sigma)$ of dimension 1. Every $W_\sigma$ is the kernel of some homomorphism $\xi_\sigma
\in \hom(\sfh^0(Y, K \otimes Q_\sigma), \cc) \cong \sf{Ext}^1(Y; Q_\sigma, \scro)$ defined
uniquely up to a non-zero factor, so that the associated extension $0 \to \scro
\to E_\sigma \to Q_\sigma \to 0$ is well-defined. Setting $X_\sigma \deff \pp(E_\sigma)$ we obtain a
holomorphic family of minimal ruled surface which contains any generic $X$.
Moreover, on every $X_\sigma$ we obtain the curves $C_{0,\sigma}, \ldots, C_{3,\sigma}$
corresponding to the subbundle $\scro \subset E_\sigma$ and the components $\{y_1, y_2,
y_3\}$ of $\sigma$. The monodromy along the family $\{X_\sigma\}_{\sigma \in S}$ acts as the full
symmetric group $\sym_3$ on the curves $C_{1,\sigma}, C_{2,\sigma}, C_{3,\sigma}$ since such
is the monodromy action on the points $y_1, y_2, y_3$. To obtain further
permutations of $C_{0,\sigma}, \ldots, C_{3,\sigma}$ we interchange $C_{0,\sigma}$ with one of the
remaining $C_{i,\sigma}$, $i=1,2,3$. For this purpose we take the line subbundle
$L_{1,\sigma} \subset E_\sigma$ corresponding to the curve $C_{1,\sigma}$ and consider the
extension $0 \to \scro \to E_\sigma \otimes L_{1,\sigma}\inv \to Q^{(1)}_\sigma \to 0$ where $Q^{(1)}_\sigma
\deff E_\sigma \otimes L_{1,\sigma}\inv / \scro$ is the quotient bundle. The same deformation
construction as above can be applied to the new extension. Consequently, the
whole monodromy action on the curves $C_{0,\sigma}, C_{1,\sigma}, \ldots,C_{3,\sigma}$ is the full
symmetric group $\sym_4$.

\medskip 
Now consider the case $g=3$. Denote by $\xi$ the element of $\sf{Ext}^1 (Y; Q,
\scro) \cong \sfh^1(Y, Q\inv) \cong \big( \sfh^0(Y, K \otimes Q) \big)^*$ defining the
extension \eqqref(2.5).  It was shown in the proof of \lemma{lem2.0} that
every a line subbundle $L \subset E$ of degree $0$ different from $\scro \subset E$
corresponds to a line $\ell$ in $\pp^3 \deff \pp\big( (\sfh^0(Y, K \otimes Q) )^*
\big)$ passing though $2$ points on $\phi_{K \otimes Q}(Y)$ and the point $[\xi]$. The
linear projection $\pi_{[\xi]}: \pp^3 \bs \{[\xi]\} \to \pp^2$ from $[\xi]$ establishes
the 1-1-correspondence between such lines $\ell$ and double points on the sextic
($d\deff \deg(K \otimes Q) =3g-3=6$) $\pi_{[\xi]}(\phi_{K \otimes Q}(Y)) $of genus $g=3$ in
$\pp^2$. The genus formula insures the existence of $\delta = \frac{(6-1)(6-2)}{2}
 -3 = 7$ nodal points. This corresponds to $7+1 =8$ curves of degree $1$ and
self-intersection $[C]^2 = 2$ on a generic minimal ruled surface $X$ over a
curve $Y$ of genus $3$. We contend that the monodromy along an appropriate
family of deformations of $X$ acts as the full symmetric group $\sym_8$ of
permutations of such curves. To show this, let us first observe that in the
case $g\geq3$ the curve $\phi_\xi(Y) \subset \pp^{2g-4}$ allows to restore $Y$ and the whole
extension \eqqref(2.5).  Indeed, the base curve $Y$ is simply the
normalization of the image, so that $\phi_\xi: Y \to \pp^{2g-4}$ can be considered as
the normalization map. The bundle $Q$ is restored from the equality $K \otimes Q =
\phi_\xi^*(\scro_{\pp^{2g-4}}(1))$ and the kernel $W_\xi$ of $\xi$ as the
$\phi_\xi$-pre-image of the hyperplane linear system on $\pp^{2g-4}$ in the full
linear linear system of $K \otimes Q$. Further, we use the fact that the monodromy
along the family $\scrz_\nu(\pp^2, [dH])$ of all nodal curves of degree $d$ in
$\pp^2$ with $\nu$ nodes acts as the full symmetric group $\sym_\nu$ of
permutation of nodes. Associating the nodes of the irreducible planar sextic
$\phi_\xi(Y)$ with $7$ of $8$ curves on the ruled surface $X$ in question we obtain
the the full symmetric group $\sym_7$ of permutation of all the curves except
the curve $C_0$ corresponding to the subbundle $\scro \subset E$ in \eqqref(2.5). To
generate the whole group $\sym_8$ we interchange $C_0$ with one of the
remaining curves as it was done in the case $g=2$.

\medskip

\newsubsection[sev-hirz]{Severi problem for Hirzebruch surfaces} We start with
the lemma which allows to ``transfer'' the local results of \refsection{sec:1}
to global families.

\newlemma{lem2.5} Let $Z$ be a finite-dimensional analytic subset of
$\scrz_d(\Delta)$ and $C^* \in Z \subset \scrz_d(\Delta)$ a proper curve whose discriminant
$\dscr(P_{C^*})$ vanishes only at $z=0$. Assume that that the codimension of
$\scrz_d^\es(C^*, 0) \subset \scrz_d$ equals the codimension of $Z^*\deff
\scrz_d^\es(C^*, 0) \cap Z \subset Z$. Then every component $Y$ of $Z_\nu \deff Z \cap
\scrz_{d,\nu}$ 
\begin{enumerate}
\item the codimension of\/ $Y \subset Z$ in $\nu$;
\item a generic curve $C \in Y$ is nodal with exactly $\nu$ nodes;
\item $Y$ contains a maximal nodal deformation of $C^*$.
\end{enumerate}
\end{lem}

\proof We use the properties of BASFD's listed in \propo{prop1.3a}, especially
the property \slvp, (see also \slsf{Remark} after \propo{prop1.3a}).
Estimating codimensions, we conclude that if $\scry \subset \scrz_d$ is an analytic
subset of pure codimension $k$ which contains $\scrz^\es_d(C^*, 0)$, then the
set $\scry \cap Z$ is also of pure codimension $k$ in $Z$. The first two
assertions of the lemma are obtained if we apply this to the loci
$\scrz_{d,\nu}$ and $\scrz_{d,\nu} \bs \scrz_{d,\nu}^\circ$, respectively. 

The last claim follows similarly. First, we apply the codimension argument to
the sequence $\scry_0 \supset \scry_1 \supset \cdots \supset \scry_{k^*}$ used in \lemma{lem1.11}.
This provides that a generic curve in $Y_k \deff Z \cap \scry_k$ has properties
(1a--1d) in the cases $k=0 \ldots k^*-1$ and properties (2a--2d) in the case $k=
k^*$. Then the argument is applied to the locus $\scrw$ which appears in the
proof of \refthm{thm1.9}.
\qed

\medskip 
\proofp of \refthm{thm0.1} follows almost immediately from the lemma.  Fix the
coordinate system $(z,w)$ on $\bff_k$ in which $z=[z_0 : z_1]$ is a projective
coordinate on the base of the ruling $\pr: \bff_k \to \pp^1$ and $w= [w_0: w_1]$
indices a projective coordinate on each fiber $\pr\inv(z)$ such that the
``infinity'' section $C_\infty$ with $C_\infty^2 = -k$ is given by the equation $w_1=0$.
In these coordinates any curve $C \subset \bff_k$ is the zero of the Weierstra\3
polynomial $P_C= \sum_{i=0} ^d a_i(z) w_0^{d-i} w_1^i$ where $a_i \in \sfh^0(
\pp^1, \scro(f + ik)$. It is easy to show that in this case $C$ belongs to the
linear system of $d \cdot C_0 + f \cdot F$ where $F$ is a fiber of the projection
$\pr: \bff_k \to \pp^1$ and $C_0$ is a section with $C_0^2 =k$. We treat
elements of $\sfh^0( \pp^1, \scro(l))$ as polynomials $a(z)$ of degree at most
$l$.

Now let $C \subset \bff_k$ be a nodal curve without multiple components. Then there
exists a fiber $F_{z^{\dag}} \deff \pr\inv(z^{\dag})$ which meets $C$ transversally.
Changing the coordinate $z$, if needed, we may assume that $z^{\dag} =\infty$. Consider
the toric action of $\cc^*$ on $\bff_k$ given by the formula $(\lambda; z, w)
\overset{\vphi} \longmapsto (\lambda\,z, \lambda^kw)$. The lift of this action on Weierstra\3
polynomials $P_C= \sum_{i=0} ^d a_i(z) w_0^{d-i} w_1^i$ with $a_i \in \sfh^0(
\pp^1, \scro(f + ik)$ is given by 
\[
P_C= \sum_{i=0} ^d a_i(z) w_0^{d-i} w_1^i \overset{\vphi_\lambda} \longmapsto 
\sum_{i=0} ^d \lambda^{-(f+ik)} a_i(\lambda\,z) w_0^{d-i} w_1^i.
\]
It is obvious that $\vphi_\lambda(C)$ are nodal for all $\lambda \neq 0$ and that the limit
curve $C^* \deff \lim_{\lambda \to 0} \vphi_\lambda(C)$ consists of the fiber $F_0 \deff
\pr\inv(0)$ with some multiplicity $m_0^*$ and $d$ sections $C_1^*, \ldots, C_d^*$
given by the equations $w = \alpha_i z^k$ with pairwise distinct $\alpha_i\in \pp^1$.
Moreover, the parameters $\alpha_i$ are exactly the $w$-coordinate of the
intersection points of $C$ with the fiber $F_\infty$. In the case $a_j = \infty$ the
curve $C^*_j$ is the ``infinity'' section $C_\infty$. There exists at most one such
section $C^*_j$. The $m_0^*$ can be computed  via the intersection index of $C$ with
$C_\infty$. It is easy to show that $m_0^* = f +k$ if $C_\infty$ is a component of $C^*$
and $m_0^* = f$ otherwise.

Now, the claim of the theorem is a special case of \lemma{lem2.5} with $Z =
\scrm( \bff_k, d,f, g=0)$ and the curve $C^* = \lim_{\lambda \to 0} \vphi_\lambda(C)$.
\qed

\medskip 
\proofp of \slsf{Theorem 1.} We maintain the notation introduced in the proof
of \refthm{thm0.1}. Let $C$ be an irreducible nodal curve on $\bff_k$. Then
$C$ lies in the linear equivalency class of $d \cdot C_0 + f \cdot F$. If $C_0$
differs from $C_\infty$, then $f = C \cdot C_\infty \geq 0$. Thus the fiber $F_\infty$ in the proof
of \refthm{thm0.1} can be chosen in that way that $C \cap F_\infty$ is disjoint from
$C_\infty$. This means that $C_\infty$ is not a component of $C^*$.

A maximal nodal perturbation $C^\times$ of $C^*$ consists of $f$ vertical lines
$F_{z_i^\times} = \pr\inv(z^\times_i)$, $i=1, \ldots ,f$, and $d$ sections $C^\times_j$, $j=1, \ldots
,d$ which meet transversally at pairwise distinct points. Thus we obtain
$k\frac{d(d-1)}{2} + d\cdot f$ nodes on $C^\times$. To denote these nodes we set $\{
x^\times_{ij} \} = F_{z_i^\times} \cap C^\times_j$ and $\{ x^\times_{ij1}, \ldots , x^\times_{ijk} \} = C^\times_i \cap
C^\times_j$.

By \refthm{thm0.1}, every component of the variety $\scrm^\circ(\bff_k, d, f, g)$ can
be reached by smoothing an appropriate subset $\bfs$ of the whole collection
$\bfn^\times \deff \Sing(C^\times)$ of the nodes of $C^\times$. Observe that this can be done
in two steps. First, we smooth an appropriate subcollection $\bft \subset \bfs$ of
$d+f-1$ nodes such that obtained curve $C^{\dag}$ is rational and irreducible, and
then smooth an appropriate collection of $g$ nodes on $C^{\dag}$. Consequently, the
theorem follows from the following lemma.

\newlemma{lem2.6} \sli The variety $\scrm^\circ(\bff_k, d, f, 0)$ is irreducible.

\slii The monodromy along $\scrm^\circ(\bff_k, d, f, 0)$ acts on the set of nodes
of a given curve $C^{\dag} \in \scrm^\circ(\bff_k, d, f, 0)$ as the full symmetric group.
\end{lem}

\proof We still maintain the notation introduced above. Furthermore, we make
geometric constructions in the affine chart of $\bff_k$ corresponding to the
finite values of the coordinates, $(z,w) \in \cc^2$. This means that a curve $C$
in the linear system of $C^*$ is defined by a usual polynomial $P_C(z,w)=
\sum_{i=0} ^d a_i(z) w^i$ of the affine coordinates $(z,w) \in \cc^2$ of the form
$P_C(z,w)= \sum_{i=0} ^d a_i(z) w^{d-i}$ such that each $a_i(z)$ is a polynomial
of degree $\leq f + ki$. In particular, each component $C^\times_j$ is given by the
equation $w= p_j(z)$ with a polynomial $p_j(z)$ of degree $\leq k$. The nodal
points $\{ x^\times_{ij1}, \ldots , x^\times_{ijk} \} = C^\times_i \cap C^\times_j$ correspond to the zeros
of the polynomial $p_i(z) -p_j(z)$. This shows that the monodromy group $G^\times$
along the locus of maximal nodal deformation contains the
$\frac{d(d-1)}2\,$-fold product of the symmetric groups $\sym_k$ of permutations
of the sets $\{ x^\times_{ij1}, \ldots , x^\times_{ijk} \}$, $1\leq i< j\leq d$. It is also obvious
that the natural action of $G^\times$ on the sets $\{C^\times_1, \ldots , C^\times_d\}$ and $\{
F_{z^\times_1}, \ldots , F_{z^\times_f} \}$ defines an epimorphism onto the product $\sym_g \times
\sym_f$ whose kernel is exactly the product $\big( \sym_k \big){}^{d(d-1)/2}$ above.
This gives a complete description of the action of the monodromy group on the
set of nodes of $C^\times$.

Now consider the action of the monodromy group on the set of nodes of some
fixed curve $C^{\dag} \in \scrm^\circ(\bff_k, d, f, 0)$. As it was shown above, we may
assume that $C^{\dag}$ is obtained from $C^\times$ by smoothing an appropriate subset
$\bft$ of the whole collection $\bfn^\times$ of nodes of $C^\times$. By monodromy
argument, we can identify $\bfn^\times \bs \bft$ with the nodes of $C^{\dag}$. 

Assume that $k\geq1$. Fix three components $D_1, D_2, D_3$ of $C^\times$ and three
nodal points $q_1, q_2, q_3 \in \bfn^\times$ such that $q_1 \in D_2 \cap D_3$, $q_2 \in D_1
\cap D_3$, $q_3 \in D_1 \cap D_2$. In particular, at most one of the components $D_1,
D_2, D_3$ is a vertical line $F_{z^\times_i}$. Observe that at most two of the
points $q_1, q_2, q_3$ belong to $\bft$ since otherwise $C^{\dag}$ would be not
rational. Assume that $q_1 \in \bft$ and denote by $C^{\ddag}$ the curve obtained from
$C^\times$ by smoothing $q_1$.  As above, we identify the set of nodes of $C^{\ddag}$
with the set $\bfn^\times \bs \{ q_1\}$. We contend that the monodromy group along
the variety of equisingular deformations of $C^{\ddag}$ contains the transposition
of $q_2$ and $q_3$. To show this, we first bring the points into $q_1, q_2,
q_3$ a ordinary triple point, then smooth the point $q_1$ creating a nodal
irreducible rational curve $D_{23}$ from $D_2$ and $D_3$, and finally move
$D_3$ in such a way that $q_2$ and $q_3$ collapse into a tangency point of
$D_1$ and $D_{23}$.  Then the monodromy around the locus of equisingular
deformations of such a constellation with tangency is the desired
transposition of $q_2$ and $q_3$. It is convenient to control the creation of
such a constellation in the affine chart introduced above. Let us call a
constellation of $D_1, D_2, D_3$ and points $q_1, q_2, q_3$ a \slsf{triangle}
and the constructed operation the \slsf{transposition of $q_2$ and $q_3$ with
 support in $q_1$}.

We contend that applying this operation with appropriate constellations of
$q_1, q_2, q_3$ we can obtain a new subset $\wt\bft \subset \bfn^\times$ with the
following properties: 
\begin{itemize}
\item[(T1)] all nodal points in $\wt\bft$ lie on the component $C^\times_1$;
\item[(T2)] smoothing the nodes in $\wt\bft$ gives a curve $\wt C^{\dag}$ lying in
 the same component $\scrm^\circ(\bff_k, d, f, 0)$ as $C^{\dag}$.
\end{itemize}
To assure the second property, it is sufficient to reach $\wt\bft$ by a chain
$\bft \ddef \bft_0 \to \bft_1 \to \bft_2 \to \cdots \to \bft_n \deff \wt\bft$ such that at
each step $\bft_i \to \bft_{i+1}$ the set $\bft_{i+1}$ is obtained from $\bft_i$
by applying a transposition supported in some point of $\bft_i$. Now assume
that at some stage we have obtained a collection $\bft_i \subset \bfn^\times$ with
property (T2) for which property (T1) fails. Then (T2) insures that there
exist a component of $C^\times$, say $D_1$, and two nodal points on $D_1$, say
$q_2$ and $q_3$, such that both $q_2, q_3$ belong to $\bft_i$ and $q_3 \in
C^\times_1$. Set $D_2 \deff C^\times_1$ and denote by $D_3$ the component of $C^\times$
passing through $q_2$. Let $q_1$ be a point from $D_2 \cap D_3$. The intersection
of $D_2 = C^\times_1$ with $D_3$ can not be empty since $k \geq 1$. Thus we obtain a
triangle constellation $D_1, D_2, D_3;\; q_1, q_2, q_3$ such that $q_1 \not\in
\bft_i$, since otherwise smoothing of $C^\times$ in $q_1, q_2, q_3$ would give a
curve of genus $1$. Performing the transposition of $q_1, q_2$ with support in
$q_3$ we obtain a new collection $\bft_{i+1}$ which have more points on
$C^\times_1$ than $\bft_i$. 

\smallskip 
To finish the proof of the lemma in the case $k\geq1$ being considered, it
remains to observe that all possible transpositions of points $q_2, q_3 \not\in
\wt\bft$ supported in $q_1 \in \wt\bft$ combined with the subgroup of the group
$G^\times$ leaving the set $\wt\bft$ invariant generate the full symmetric group of
permutation of the set $\bfn^\times \bs \wt\bft$. The details are left to the
reader.

\medskip
Finally, consider the case $k=0$. Since $\bff_0 = \pp^1 \times \pp^1$, the
projections on the first and on the second factor are two possible rulings,
$\pr: \bff_0 \to \pp^1$ and $\pr': \bff_0 \to \pp^1$. The curve $C^\times$ consists of
$d$ vertical fibers $F_i \deff \pr\inv(z^\times_i)$ and $f$ horizontal fibers $F'_j
\deff \pr^{\prime\;-1}(w^\times_j)$. This gives $d\cdot f$ nodal points $x^\times_{ij} \deff F_i
\cap F'_j$. The nomodromy group $G^\times$ along the family of equisingular
deformations of $C^\times$ permutes the horizontal and vertical fibers
independently, so $G^\times = \sym_d \times \sym_f$. Our working procedure in the case
$k=0$ involves a \slsf{rectangular constellation}, consisting of two vertical
fibers $D_1, D_2$, two horizontal fibers $D'_1, D'_2$, and the edges $q_{ij}
\deff D_i \cap D'_j$. For such a constellation $D_i, D'_j, q_{ij}$, we smooth two
points lying on one of the sides of the rectangular, say $q_{11}$ and $q_{12}$,
creating an irreducible rational curve $D$ in the linear system $D_1 + D'_1 +
D'_2$. Moving the remaining side $D_2$ around the locus of tangency of $D_2$
with $D$ we obtain a family whose monodromy transposes the points $q_{21}$ and
$q_{22}$. We call it the \slsf{transposition of $q_{21}$ and $q_{22}$ with
 support in $q_{11}$ and $q_{12}$}. We contend that there exists a chain 
$\bft \ddef \bft_0 \to \bft_1 \to \bft_2 \to \cdots \to \bft_n \deff \wt\bft$ of subsets 
of the set $\bfn^\times$, such that each $\bft_{i+1}$ is obtained from $\bft_i$ by a
transposition with support in $\bft_i$, and such that the final collection
$\wt\bft$ has the following properties:
\begin{itemize}
\item[(T1')] all nodal points in $\wt\bft$ lie on the components $F_1$ and $F'_1$;
\item[(T2)] smoothing the nodes in $\wt\bft$ gives a curve $\wt C^{\dag}$ lying in
 the same component $\scrm^\circ(\bff_k, d, f, 0)$ as $C^{\dag}$.
\end{itemize}
As above, this implies the irreducibility of $\scrm^\circ(\bff_k, d, f, 0)$.
Details of the proof of the existence of such a chain $\bft = \bft_0 \to \bft_1 \to
\bft_2 \to \cdots \to \bft_n = \wt\bft$ and of the second part of the lemma in the case
$k=0$ are left to the reader.
\qed

\bigskip


\ifx\undefined\bysame
\newcommand{\bysame}{\leavevmode\hbox to3em{\hrulefill}\,}
\fi

\def\entry#1#2#3#4\par{\bibitem[#1]{#1}
{\textsc{#2 }}{\sl{#3} }#4\par\vskip2pt}

\def\noentry#1#2#3#4\par{}

\def\mathrev#1{{{\bf Math.\ Rev.:\,}{#1}}}


\end{document}